\documentstyle[11pt,amsmath,amsfonts,here ]{article}

\newtheorem{th}{Theorem}

\def\R{{\mathbb R}}
\def\C{{\mathbb C}}
\def\Z{{\mathbb Z}}

\def\a{\alpha}
\def\b{\beta}
\def\e{\varepsilon }
\def\g{\gamma}
\def\o{\omega}
\def\i{\iota}
\def\k{\kappa}

\def\t{\tau}
\def\rw{\rightarrow }

\def\d{\delta}

\def\l{\lambda}
\def\s{\sigma}
\def\n{\newline}
\def\w{\newline \newline}
\def\O{{\cal O}}
\def\M{{\cal M}}
\def\H{{\cal H}}
\def\G{{\cal G}}

\begin{document}
\title{\large \bf 
UNIFORM ASYMPTOTIC BOUND ON THE NUMBER OF ZEROS OF ABELIAN INTEGRALS}
  \author{\large Alexei Grigoriev\footnote{Sector of Functional Analysis, 
SISSA, Trieste. Email: alexg@sissa.it.}}
  \date{}
  \maketitle
  \begin{abstract}
We give a uniform asymptotic bound for the number of zeros of complete 
Abelian integrals in domains bounded away from infinity and 
the singularities. 
\end{abstract}
  
\bigskip
\bigskip
\bigskip
{\bf 0. Introduction.}\n

In [Aea], V.I. Arnol'd posed the question of the number of 
limit cycles which can arise 
by perturbing an integrable planar vector field of degree $d$ by a small
polynomial perturbation of the same degree. 
He indicates the 
relation of this question to the determination of the number of
zeros of complete Abelian integrals parameterized by a value of a 
polynomial of degree not greater than $d$, which is just the Hamiltonian
in the case that the vector field is a Hamiltonian vector field (which is the 
case treated in the present article). 
It is pointed out that in this case the number of zeros of these 
integrals (in an interval which does not contain a critical value of the 
Hamiltonian, for example), is always finite and could be bounded by a number 
which depends only on the degree $d$
(this key result is due to Khovanskii and Varchenko [Kh], [Var]).
It is noted that no effective estimates of this bound are known.\n

The main result of this article is formulated in Theorem 0.2 below
(proved as Corollary 4.5 and Corollary 4.7 in the text). It
was initially obtained in the author Ph.D. thesis [Gr], and
gives a (doubly exponential in $d$) 
bound on the number of zeros of complete Abelian integrals parameterized by the
value of a Hamiltonian, in an interval whose distance from the 
critical values of the Hamiltonian is fixed, and which is contained
in a fixed bounded interval. The bound is uniform in the sense that
it holds for all Hamiltonians in a dense open
subset of the space of polynomials of degree $d$. In fact, the bound 
is stated for the number of zeros of the integrals in domains in $\C$
which satisfy certain natural restrictions, and which are likewise 
bounded away from the critical values of the Hamiltonian and the 
infinity. This uniformity of the asymptotic bound (at present in domains 
bounded away from the singularities and infinity) appears to overcome
the corresponding difficulty in the methods suggested in 
[NY1], [GI].\n 
 
We now give a more precise description of the problem. 
Consider a polynomial $H(x,y)\in\R[x,y]$ of degree $d$. Suppose that
for some $a,b\in\R$, $(a,b)$ does not contain critical values of $H$,
and let $\g_{t_0}$ denote a compact component of $H^{-1}(t_0)$ for some
$t_0\in(a,b)$. For each $t_1\in[a,b]$, 
let $\g_{t_1}$ denote the compact component of $H^{-1}(t_1)$ which is 
obtained from $\g_{t_0}$ by the obvious 
isotopy as $t$ varies from $t_0$ to $t_1$. Let 
$\s(t):(a,b)\rw \R^2$ denote an analytically embedded
segment, such that $\s(t)\in\g_t$ $\forall t\in(a,b)$, and $\s$ is 
transversal to $\g_t$
for all $t\in(a,b)$. Let $P(x,y),Q(x,y)$ be polynomials of degree not greater
than $d$, and consider the following $\e$-dependent  perturbation of the
Hamiltonian vector field determined by $H$:
$$
0.1) \ \ \ \ 
\dot x={\partial H\over\partial y}+\e Q, \ \ \ \ \dot y=-{\partial H
\over \partial x} -\e P.
$$
For $\e=0$ every trajectory which starts at $\s(t)$, $t\in(a,b)$, 
returns to $\s(t)$. For any $[a',b']\in(a,b)$, there exists a 
small enough nonzero $\e_0>0$, such that 
the map from $[-\e_0,\e_0]\times \s[a',b']$ to $\s(a,b)$, defined
for any $-\e_0\leq \e\leq\e_0$ as the corresponding Poincare return map
from $\s[a',b']$ to $\s(a,b)$, indeed exists, and is moreover analytic.
Fix $\e$, $-\e_0\leq\e\leq\e_0$; for any $t\in [a',b']$ suppose
that the Poincare map sends $\s(t)$ to $\s(t')$. Then the {\em displacement
map} $\delta(\e,t)$ is equal by definition to $t-t'$ (so that 
$\delta(0,t)\equiv 0$ as a function of $t$). So  $\delta(\e,t)$ is an analytic
function on  $[-\e_0,\e_0]\times [a',b']$. Observe that for any fixed
$-\e_0\leq \e\leq\e_0$, the zeros of $\delta(\e,\cdot)$ correspond to
periodic trajectories of the perturbed vector field, crossing $\s[a',b']$.\n

It is well known (cf. [R], pg. 72-73, for example) that the $\e$-derivative of $\d(\e,t)$
on the line $\e=0$ is given by
$$
0.2) \ \ \ \ -\int_{\g_t}P(x,y)dx+Q(x,y)dy
$$ 
(the orientation of $\g_t$ being given by the flow of the Hamiltonian 
vector field).
Observe that as $\e$ varies, starting from $\e=0$, zeros of $\delta(\e,\cdot)$
may bifurcate from $\{0\}\times\s[a',b']\subset\R^2$ at point $(0,\s(t))$
only if $(0,t)$ is a nonsmooth point of the (semianalytic) set 
$\{(\e,t)\in [-\e_0,\e_0]\times[a',b']: \delta(\e,t)=0\}$. 
Such points must be the critical points of
the map $\d(t,\e)$ on the line $\e=0$; since the derivative of 
$\d(t,\e)$ on the line $\e=0$ is given by 0.2), $(0,t)$ is such a critical 
point iff 0.2) is zero. Therefore the number of bifurcation points
is bounded by the number of zeros of the integral 0.2) on $(a,b)$.

In other words, if 0.2) is not identically zero as a function of $t$ and
its number of zeros on $(a,b)$ is $N$, for any $[a',b']\subset(a,b)$ there
exists $\e_0>0$, such that for any $\e$, $0<|\e|\leq|\e_0|$, the number of
limit cycles (i.e. isolated periodic trajectories) of the vector field 
0.1) which cross $\s[a',b']$ is bounded by $N$.\n

Since the time [Aea] was published, numerous partial results became available
regarding effective estimates of the number of zeros of Abelian integrals,
especially in low degrees; we refer the reader to the survey
in [I2]. For the general case, however, there still exist 
no effective estimates. Below, we formulate the main result 
of this article. 
Since we consider the Abelian integrals in the complex plane instead 
of the real line, we give first some additional background.

For any polynomial $H(x,y)\in\C[x,y]$, 
there exists a finite set $\Sigma_H\subset\C$,
whose points are called the atypical points of $H(x,y)$, 
such that $H:\C^2-H^{-1}(\Sigma_H)\rw\C-\Sigma_H$ is a (smooth, not 
holomorphic in general) locally trivial
(fiber) bundle. 
The integral 0.2) can be then analytically continued 
along any path not passing through the atypical points (which are 
generically just the critical values of $H(x,y)$). It is therefore
a (transcendental) multivalued analytic function on $\C-\Sigma_H$.
This fact generalizes as follows. Any element of the homology group of
a fiber over a point which is not atypical, can be naturally continued 
along any path not passing through the atypical points, the resulting
homology class being dependent 
only on the homotopy class of the path connecting 
the initial and final points. The result of such
(multivalued) continuation is called
a continuously varying cycle, since in the case that 
the initial element of the homology group
is realized by a map of $S^1$ into the fiber over the initial point,
the continuation along a path not passing through the atypical points
is realized by any isotopy of $S^1$ which is contained in the preimage 
of the path and respects the fibers. The 
integral of $Pdx+Qdy$ over any continuously varying cycle 
is a multivalued analytic function.\n 

Below, $\H^d\cong\C^{(d+1)(d+2)/2}$ denotes the subspace of $\C[x,y]$ 
consisting of polynomials of degree not greater than $d$.
The bound is given for polynomials $H(x,y)$ whose highest homogeneous 
part is a product of pairwise different (up to multiplication by a nonzero
constant) linear factors. Equivalently, these are polynomials of degree $d$,
whose level curves intersect the line at infinity at precisely $d$ points.
We say in this case that the polynomial $H(x,y)$ is regular at infinity.
For such polynomials, the only atypical points are their critical values.\n

{\bf Definition 0.1.} Let $\{t_1,..,t_N\}\subset\C$ be a finite set of points.
A domain $U\subset\C- \{t_1,..,t_N\}$ is called a 
{\em simple domain in $\C- \{t_1,..,t_N\}$} if 
there exist $N$ nonintersecting rays $r_1,..,r_N$, issuing from the points  
$t_1,..,t_N$ respectively, such that $U\subset\C-\cup_i r_i$.\n

{\bf Theorem 0.2 [Gr].} {\it There exists a universal constant $c>0$,
such that the following holds.

Let $H\in\H^d$ be regular at infinity. 
Let $P,Q\in\H^d$. Let $\rho>0$ be any positive number, and
let $U$ be a simple domain in $\C-\Sigma_H$, contained in the unit disc,
whose distance to $\Sigma_H$ is at least $\rho$. Then the number of zeros
in $U$
of the integral $\int_{\g(t)}P(x,y)dx+Q(x,y)dy$, where $\g(t)$ 
is any continuously varying 
cycle in the locally trivial bundle determined by $H$, is not greater than }
$$
\left({2\over \rho}\right)^{2^{d^c}}.
$$  

\bigskip

{\bf Remark 0.3.} The {\em existence} of a bound 
for the number of zeros in {\em any}
simple domain in $\C-\Sigma_H$, follows from Theorem 
5 in [Kh].\n  

{\bf Remark 0.4.} Since the title of [IY] mentions a double exponential
estimate as well, the following comment might be helpful.
In [IY], the polynomial $H(x,y)$, of some degree $d_0$, 
is {\em fixed}, the polynomials $P,Q$ being
any polynomials of degree $d$. In this setting, it is now known
(Khovanskii and Petrov, yet unpublished) that the bound is in fact {\em linear}
in $d$.\n 

We now outline the contents of the article. In section 1 we deal with 
linear differential equations whose
coefficients are meromorphic and depend on parameters, 
and which are not singularly perturbed along any holomorphic arc
in the parameter space. We say that such equations depend {\em regularly} 
on parameters (this terminology is due to the author; we have not found
an analogous notion in the literature). 
It is shown that solutions of such equations admit
a uniform (over compact subsets in the parameter space) bound for
the number of zeros in domains of the type considered in
Theorem 0.5. I thank S. Yakovenko for the suggestion
to use his theorem ([Y]) to simplify some of the arguments. 
To exclude a linguistic confusion,
we remark that the notion of regular dependence on parameters  
and the notion of a polynomial regular at infinity are {\it not} related.  

In section 2, we point out that if the meromorphic coefficients of a linear 
differential equation are
quotients of {\em integral} polynomials of known degree and height, the bound
of section 1 can be made effective. We rely on a theorem of Renegar
which gives effective estimates for elimination of quantifiers in the 
first order theory of the reals.

In section 3, we construct such equations for the Abelian integrals, relying
on a result of Gavrilov ([Ga],[N]) and a quantitative version of a theorem of
Ilyashenko [I1], the parameters being the coefficients of
the polynomial $H$ and (roughly speaking) the coefficients 
of the polynomials $P$ and $Q$. In section 4, we show that the linear
differential equations constructed indeed depend regularly on parameters.
Theorem 0.2 then follows from the result of
section 2 and from an additional argument, which shows that Picard-Fuchs
equations stay regularly dependent on parameters also after 
algebraic parameter dependent changes of variable.\n

{\bf Remark 0.5.} One may show using similar arguments [Gr], that the number
of zeros in the unit disc of (components of) solutions of the linear
differential system $\dot x=(A_0+A_1 t+...+A_d t^d)x$, where $A_0,..,A_d$
are some $n\times n$ matrices over $\C$, of (say $l_\infty$) 
norm not greater than 1, is bounded by $2^{2^{(nd)^{c}}}$, where 
$c>0$ is again some universal constant.
This improves the bound which follows for this case from [NY2], Theorem 1.\n 

This article is in essence a presentation of the major part of 
the results of [Gr]. We have tried to make the exposition here
clearer. In particular, some of the arguments in the thesis are made simpler.
The construction of the Picard-Fuchs system 3.15) for the Abelian integrals,
while being based on the same elementary idea, 
differs somewhat from the one given in the thesis, 
where there was a gap in the 
proof of the corresponding theorem. Recently, there appeared  
another such construction [N] (see also [NY3]), 
giving a much better information 
on the system than that which is provided by our argument. 
It seems, however, that using this construction here 
would not improve the doubly exponential bound of Theorem 0.2.\n

Finally, an important remark about notations. Instead of introducing each time
a new universal constant when such is needed, we write $O(x)$ to denote
$c_1 x+c_2$, where $c_1,c_2>0$ are some {\em universal} constants (here $x$ is
any positive number). Thus, $O(1)$ simply denotes a universal constant.\w

{\bf Acknowledgments.} I thank my thesis adviser, S. Yakovenko, for 
introducing me to the problem treated in the text, and him and D. Novikov for 
stimulating discussions. I thank A. Glutsuk, Y. Il'yashenko and  
the anonymous referees of the thesis for their remarks.\w

{\bf 1. Linear differential equations which depend regularly on parameters.}

\bigskip

Below, we formulate both the definitions and the statements for polydiscs  
instead of more general domains, as no greater generality is needed.
We denote by ${\cal O}(T)$ the set of holomorphic functions on the set
(open or closed) $T$, and by ${\cal M}(T)$ the set of meromorphic functions
on $T$.\n

Let $U\times W\subset \C\times\C^p$, $p\geq 1$, be a polydisc, and let the 
coordinates on $\C\times\C^p$ be denoted by $(t,\l_1,..,\l_p)=(t,\l)$. Below
we think of $\l_1,..,\l_p$ as parameters. Let $f(t,\l)\in{\cal M}(U\times W)$
be a meromorphic function; we consider it as a parameter dependent 
meromorphic function of $t$. To emphasize its domain of definition, we
will sometimes say that $f(t,\l)$ is a parameter dependent 
meromorphic function on $U\times W$.

We say that $f(t,\l)$ depends {\em holomorphically} on $\l$ at 
$\l_0=(\l_{10},..,\l_{p0})\in W$, if $f(t,\l)=a(t,\l)/b(t,\l)$ where
$a(t,\l),b(t,\l)\in{\cal O}(U\times W)$, such that $b(t,\l_0)$ is not
identically zero (it is suitable to recall here that a meromorphic function
on a polydisc is always a quotient of two holomorphic functions).\n

{\bf Definition 1.1.} A {\em holomorphic arc} in the parameter space,
passing through $\l_0\in W$, is a germ of a holomorphic map from
a neighbourhood of $0\in\C$ into $\C^p$, mapping $0$ to $\l_0$.\n

{\bf Definition 1.2.} A parameter dependent meromorphic function $f(t,\l)$
on a polydisc 
$U\times W\subset \C\times\C^p$ will be said to depend {\em regularly} on 
parameters at $\l=\l_0\in W$, if $f(t,\l(\e))$ depends holomorphically
on $\e$ at $\e=0$, for every holomorphic arc $\l(\e)$ passing through $\l_0$,
on which $f(t,\l)$ is defined. We will say that $f(t,\l)$ depends regularly
on $\l$ in $W$, if it depends regularly on $\l$ at any point of $W$.\n

{\bf Remark 1.3.} We say that $f(t,\l)$ is defined on $\l(\e)$ if it is
possible to write $f(t,\l)$ as a quotient $a(t,\l)/b(t,\l)$, 
$a(t,\l),b(t,\l)\in {\cal O}(U\times W)$, where $b(t,\l(\e))\not\equiv 0$.\n

{\bf Remark 1.4.} As it was noted, in the case $dim \l=1$, the notions 
of holomorphic dependence and regular dependence at a fixed point of
the parameter space coincide. This is not the case for $dim \l>1$.
Indeed, consider $\l_1\l_2/(\l_1^2+\l_2^2 t)$, which depends regularly on $\l$
at $\l=0$, but not holomorphically. In section 4 we verify that the 
coefficients of Picard-Fuchs equations constructed in section 3 depend
regularly on parameters; due to complexity of computations, 
we do not have a counterexample which would 
show that at the same time they do not depend holomorphically on parameters,
though it is likely that they do not (when the dimension of the parameter
space is greater than one). In another related situation 
linear differential equations which depend regularly, but
not holomorphically, on parameters, do arise naturally. Such are, in general,
the linear differential equations for the
components of a linear differential system with a polynomial system matrix,
with parameters being the coefficients of the polynomial entries ([Gr]).\n  

{\bf Remark 1.5.} Definition 1.2, put slightly differently, is:
for {\em any} representation of $f(t,\l)$ as $a(t,\l)/b(t,\l)$, $a(t,\l),b(t,\l)\in {\cal O}(U\times W)$, and for any holomorphic arc $\l=\l(\e)$ passing 
through $\l_0$ such that $b(t,\l(\e))\not\equiv 0$, $a(t,\l(e))/b(t,\l(\e))$
depends holomorphically on $\e$. It is in fact equivalent to the following
alternative definition:
for {\em some} representation of $f(t,\l)$ as $a(t,\l)/b(t,\l)$, $a(t,\l),b(t,\l)\in {\cal O}(U\times W)$, 
and for any holomorphic arc $\l=\l(\e)$ passing 
through $\l_0$ such that $b(t,\l(\e))\not\equiv 0$, $a(t,\l(e))/b(t,\l(\e))$
depends holomorphically on $\e$.
The equivalence of both definitions 
is a corollary of Lemma 1.6 below. 
\n

{\bf Lemma 1.6.} {\it Let $f(t,\l)$ be a parameter dependent meromorphic 
function on a polydisc $U\times W\subset \C\times\C^p$, and 
suppose that for the 
holomorphic arc $\l=\l(\e)$, $f(t,\l(\e))$ is defined and does not depend
holomorphically on $\e$ at $\e=0$. Then there exists $m>0$, such that
for any holomorphic arc $\l={\widehat \l}(\e)$, 
for which $\l(\e)-{\widehat \l}(\e)=o(\e^m)$,
$f(t,{\widehat \l}(\e))$ does not depend holomorphically on $\e$ at $\e=0$ as well.}

{\bf Proof.} By assumption 
$f(t,\l(\e))$ is defined and does not depend holomorphically
on $\e$ at $\e=0$. Then 
there exists a representation $f(t,\l)=a(t,\l)/b(t,\l)$, 
$a(t,\l),b(t,\l)\in {\cal O}(U\times W)$, such that $b(t,\l(\e))\not\equiv 0$ 
and such that, writing $a(t,\l)$ and $b(t,\l)$ as the power series (convergent
in some small polydisc around $(t_0,\l_0)\in U\times W$) $a_0(\l)+a_1(\l)t+...$
and $b_0(\l)+b_1(\l)t+...$ respectively, the following holds.\n
All $b_j(\l(\e))$, $j=0,1,2,..$, have a zero of order at least $k\geq 1$ at 
$\e=0$, 
while for some $i$, $a_i(\l(\e))$ has there a zero of a smaller order 
(if at all). Clearly, for any holomorphic arc ${\widehat\l}(\e)$ for which
$\l(\e)-{\widehat \l}(\e)=o(\e^k)$, the order of vanishing of $a_i( {\widehat\l}(\e))$
and of $b_j({\widehat\l}(\e))$, $j=0,1,2,..$ stays the same. We conclude that
$f(t,{\widehat\l}(\e))$ is defined and does not depend holomorphically on $\e$
at $\e=0$.\ \ \ $\Box$\n

Lemma 1.6 will be often
used as follows. 
Suppose that
a given parameter dependent meromorphic function, restricted to a certain
holomorphic arc in the parameter space, does not depend
holomorphically on the arc parameter.  
Then there exists another holomorphic arc, such that the restriction of the 
meromorphic function to this arc again does not depend holomorphically on
the arc parameter, and which (possibly unlike the original arc) is in a 
 general
position (so that it does not lie on certain exceptional subsets of the 
parameter space $\C^p$). 

One may characterize regular dependence on parameters as follows.\n

{\bf Proposition 1.7.} {\it Let $f(t,\l)=a(t,\l)/b(t,\l)$, $a(t,\l),b(t,\l)\in {\cal O}(U\times W)$. Denote by $S\subset W$ the set $\{\l\in W: 
b(\cdot,\l)\equiv 0\}$, and let $K$ be any compact subset of $U$ with nonempty 
interior. For a fixed $\l\not\in S$, 
denote by $K(\rho,\l)$ the set obtained by removing from $K$
the points whose distance from the (discrete) set $\{t\in U:b(t,\l)=0\}$ is
smaller than $\rho$. Then the following are equivalent:\n
i) $f(t,\l)$ depends regularly on $\l$ in $W$,\n
ii) for any compact set $F\subset W$, and any $\rho>0$
$$
1.1) \ \ \ \ 
sup_{\l\in F-S} \ sup_{t\in K(\rho,\l)} \  {|a(t,\l)|\over |b(t,\l)|} \ < \ \infty
$$
whenever defined (i.e. whenever the set $\{(t,\l):\l\in F-S, t\in 
K(\rho,\l)\}$ is nonempty).}

{\bf Proof.} We prove the proposition when $a(t,\l),b(t,\l)$ are polynomials
(with complex coefficients), and then comment how it is to be modified in the
general case. So let $a(t,\l),b(t,\l)\in\C[t,\l]$, and suppose that i) holds 
but ii) does not hold.

In this case 
there exists a compact polydisc $K\times F\subset U\times W$ and a 
positive number $\rho>0$, such that 
$$
sup_{\l\in F-S} \ sup_{t\in K(\rho,\l)} \ {|a(t,\l)|\over |b(t,\l)|} \ = \ \infty.
$$   
Let $D_r\subset F$ be the set of points at (Euclidean) distance from $S$, 
$r>0$. Denote by $M_r$ the subset of $D_r$, for which
$$
1.2) \ \ \ \ sup_{t\in K(\rho,\l)} \  {|a(t,\l)|\over |b(t,\l)|} 
$$
is defined and is equal to its maximal value on $D_r$. 
If at all there are points on $D_r$ for 
which $K(\rho,\l)$ is nonempty, $M_r$ will be a nonempty set.
Denote by $\Gamma$ the set $\cup_{r>0}M_r$.

Recall now the notion of a first order formula in the 
sense of elementary mathematical logic. 
Taking the language to be the 
language of ordered rings (i.e. the language is the set of 5 symbols
$\{0,1,+,\cdot,>\}$), it is not difficult to write a 
first order 
formula defining the set $\Gamma$ (using the real and imaginary parts
of the coefficients of the given polynomials $a(t,\l),b(t,\l)$). 
By the Tarski-Seidenberg principle (by which we mean here the fact that
the first order theory of the reals eliminates quantifiers), $\Gamma$ is then 
a semialgebraic subset of $\R^{2p}\cong\C^p$. 
If $\Gamma$ is bounded away
from $S$, it means that for some $r_0>0$, $K(\rho,\l)$ is empty for all
$\l\in D_r$, $r<r_0$. 
One cannot have 
$$
sup_{\l\in F-S} \ sup_{t\in K(\rho,\l)} \ {|a(t,\l)|\over |b(t,\l)|} \ = \ \infty,
$$ 
since $|b(t,\l)|$ is bounded away from zero on $\cup_{\l\in D_r, r\geq r_0} 
K(\rho,\l)$ (and since $|a(t,\l)|$ is bounded on the compact set $K\times F$).
We conclude that $\Gamma\subset F$ has a limit point on $S\cap F$. 
By the Curve Selection lemma there exists a real analytic curve 
$\l(\e):(-\e_0,\e_0)\rw \R^{2p}\cong\C^p$, 
such that $\l((-\e_0,\e_0)-\{0\})\subset F$
and $\l(0)\in S$. This real analytic curve defines a holomorphic arc
passing through $\l(0)$, which we also denote by $\l(\e)$. It follows that 
$$
1.3) \ \ \ \ sup_{t\in K(\rho,\l(\e))}  \ {|a(t,\l(\e))|\over |b(t,\l(\e))|} 
$$
is not bounded as $\e$ tends to 0. This implies that $a(t,\l(\e))/ b(t,\l(\e))$
does not depend holomorphically on $\e$ at $\e=0$, contradicting the initial
assumption that i) is true. Therefore i) implies ii).

To show the converse, suppose that i) does not hold. Then there exists a 
holomorhic arc $\l(\e)$, for which $b(t,\l(\e))\not\equiv 0$ and
$$
{a(t,\l(\e))\over b(t,\l(\e))} \ = \ {1\over\e^s} {\a(t,\l(\e))\over \b(t,\l(\e))},
$$
for some $s\geq 1$, such that for some $t_0\in U$, $\a(t_0,\l(0))\neq 0$ and
$\b(t_0,\l(0))\neq 0$. For small enough $\rho>0$, and since $K$ has a 
nonempty interior, $K(\rho,\l)$ is nonempty for all $\l\in F-S$, where $F$ is 
a compact polydisc in $W$. Then   
$$
sup_{t\in K(\rho,\l(\e))} \ {|a(t,\l(\e))|\over |b(t,\l(\e))|} \ = \  
sup_{t\in K(\rho,\l(\e))} \ {1\over |\e|^s}{|\a(t,\l(\e))|\over |\b(t,\l(\e))|}
$$
is unbounded as $\e$ tends to 0, and therefore ii) does not hold. Thus
ii) implies i).

Regarding the changes needed to be done when $a(t,\l),b(t,\l)$ are not
necessarily polynomials, one has to replace the semialgebraic category with
the subanalytic one (cf. [BM], for example). 
The proof that the set $\Gamma\in \R^{2p}\cong\C^p$ is
subanalytic, is more cumbersome now. One may either use a 
language suitable for dealing with subanalytic sets (so that they 
(or a subclass of them) become
definable), or, what is essentially the same, 
one may replace the logical operations in the formula defining  
the semialgebraic $\Gamma$ above, by corresponding set-theoretic 
operations (being careful not to project unbounded subanalytic sets).
The other details are virtually unchanged.  \ \ \ \ $\Box$\n

We are now coming to the issue which motivated the considerations above.
Let
$$
1.4) \ \ \ \ y^{(n)}+{c_{n-1}(t,\l)}y^{(n-1)}+...+
{c_0(t,\l)}y=0
$$
be a parameter dependent 
linear differential equation with meromorphic coefficients in the 
polydisc $U\times W\subset\C\times \C^p$.\n

{\bf Definition 1.8.} The linear differential equation 1.1 will be said 
to depend regularly on parameters at $\l=\l_0\in W$ (respectively, in $W$) 
if its coefficients depend regularly on parameters at $\l_0$ (respectively,
in $W$).\n

As shown below (Theorem 1.9), 
Proposition 1.7 implies that the number of zeros
of solutions of such equations in domains satisfying some natural
restrictions, is uniformly bounded over compact sets in the parameter space.
It is the quantitative version of this assertion, formulated in section 2,
which then allows to give asymptotic bounds on the number of zeros of
Abelian integrals. 

Note, that if there exists a holomorphic arc 
$\l(\e)$, such that for some $i$, 
$0\leq i\leq n-1$ $c_i(t,\l(\e))$ is defined and does not
depend holomorphically on $\e$ (at $\e=0$), but other coefficients are
not necessarily defined on $\l(\e))$, 
then by Lemma 1.6 there exists another arc
${\widehat\l}(\e)$, ${\widehat\l}(0)=\l(0)$, 
(a perturbation of the original arc $\l(\e))$) along which all
the coefficients of 1.4) are defined, and 
$c_i(t,{\widehat\l}(\e))$ does not
depend holomorphically on $\e$. Then the equation 1.4), restricted to the arc 
${\widehat\l}(\e)$, takes the form
$$
1.5) \ \ \ \ \e^s y^{(n)}+\a_{n-1}(t,\e)y^{(n-1)}+..+\a_{0}(t,\e)y=0,
$$
with $s\geq 1$ and $\a_i(t,\e)\in{\cal O}(U'\times W')$ for all $i$,
in some polydisc $U'\times W'\subset \C\times \C$, such that
$(t_0,0)\in U'\times W'$ for some $t_0\in U$ (in fact for any $t_0$ outside
of a certain discrete subset of U, there will exist a polydisc in which
1.4) takes the form 1.5)).
Parameter dependent differential equations 
such as 1.5), are usually called {\em singularly perturbed}.
So we may rephrase Definition 1.8 as follows: we say that 1.4) depends 
regularly on $\l$ at $\l_0$, if along no holomorphic arc passing through
$\l_0$, and on which 1.4) is defined, is 1.4) singularly perturbed.\n

Let us write the equation 1.4) in the form
$$
1.6) \ \ \ \ y^{(n)}+{a_{n-1}(t,\l)\over b(t,\l) }y^{(n-1)}+...+
{a_0(t,\l)\over b(t,\l)}y=0
$$ 
where $a_i(t,\l)\in{\cal O}(U\times W)$, $i=0,..,n-1$, $b(t,\l)\in{\cal O}(U\times W)$ (again, this is always possible since $U\times W$ is a polydisc).
We denote by $S\subset W$ the set
$$
S=\{\l: b(\cdot,\l)\equiv 0\}.
$$ 
For each fixed $\l\in W-S$, we denote by $\Sigma_\l\subset\C$ the set 
of true singular points of 1.6) (i.e. the points where some solution of 1.6
is not holomorphic); $\Sigma_\l\subset Z_\l$ where $Z_\l=\{t\in U:b(t,\l)=0\}$.\n

Recall from section 0, 
that a {\em simple domain in $\C-\Sigma_\l$} is any domain $V$ contained
in a simply connected domain $G\subset \C-\Sigma_\l$ 
obtained by removing from $\C$ nonintersecting rays which 
initiate at the points of $\Sigma_\l$ and go to infinity (in other words, $G$
is obtained by making straight line cuts at the points of $\Sigma_\l$).\n

{\bf Theorem 1.9.} {\it Suppose that 1.6) depends regularly on $\l$ in $W$.
Then for any compact sets $K\subset U$, $F\subset W$, and for any $\rho>0$,
there exists $N\geq 0$, such that the following holds.

Fix $\l\in F-S$, and let $y(t)$ satisfy 1.6) for that value of the parameter.  
Then the number of zeros of $y(t)$ in any simple domain $V$ in $\C-\Sigma_\l$,
$V\subset K$, for which $dist(V,\Sigma_\l)\geq \rho$, is not greater than $N$.}
\n

In other words, the number of zeros of solutions of 1.6) in simple domains 
in $\C-\Sigma_\l$, contained in $K$ and bounded away by $\rho$ from 
$\Sigma_\l$, is uniformly bounded over $\l\in F-S$.

Theorem 1.9 is in fact the corollary of Proposition 1.7 and the following
theorem from [Y], which gives a bound on the variation of argument of a 
solution of a linear differential equation on a line segment, in terms 
of a bound on the magnitude of the coefficients of the equation on
that segment. In a sense, this is the generalization of the fact that
the number of oscillations of the solutions of $y''+K y=0$ on 
the segment $[0,1]$ is bounded by $K$ (in fact by $\sqrt{K}/2\pi$) for 
$K\geq 1$, as for $K>0$ this equation describes a linear harmonic oscillator
of frequency $\sqrt{K}/2\pi$.\n

{\bf Theorem 1.10 [Y].} {\it Let the linear differential equation 
$$
y^{(n)}+a_{n-1}(t)y^{(n-1)}+..+a_0(t)y=0
$$ 
have coefficients which are holomorphic on a finite line segment $I\subset\C$
of length $l$, their modulus being bounded there by $C\geq 1$.
Then the variation of the argument of any nontrivial solution of this equation
on $I$ is bounded by }
$$
Var \ Arg \ y(t)|_I \ \leq \ \pi(n+1)(1+lC/log(3/2)).
$$

The fact that the domain $V$ in Theorem 1.9 is simple, permits us to decompose 
it into a (uniformly bounded) number of pieces, each of which is contained
in a polygonal domain, such that the sum of the lengths of the segments 
which constitute the boundary of these polygonal domains, is uniformly
bounded, and the cardinality of the set of segments is uniformly bounded as well.
We then combine Proposition 1.7 and Theorem 1.10 to prove Theorem 1.9.
The formal argument is given below.\n

{\bf Lemma 1.11.} {\it Let $\Sigma$, $Z$ be finite subsets of a 
bounded rectangular 
domain $U\subset\C$. Let $V$ be a simple domain in $\C-\Sigma$, such that
$\overline{V}\subset U$, and such that $dist(V,\Sigma)\geq \rho$, 
$dist(V,\partial U)\geq \rho$. Then there
exist linear segments $\g_1,..,\g_N$, all contained in $U$, 
$N\leq O(|\Sigma|^2)$, such that\n
i) for all $i=1,..,N$,  $dist(\g_i,Z)\geq \rho/O(|Z|)$, 
$dist(\g_i,\partial V) \geq \rho/2$;\n
ii) let $f$ be a multivalued analytic function in $\C-\Sigma$; given a 
branch $\widetilde{f}$ on $V$, there exist branches $\widetilde{f_1},..,
\widetilde{f}_N$ defined on $\g_1,..,\g_N$, such that the number of zeros of
$\widetilde{f}$ in $V$ is not greater than the sum of absolute values of the 
variations of arguments of $\widetilde{f_1},..,\widetilde{f}_N$ on $\g_1,..,\g_N$,
respectively, divided by $2\pi$.}

{\bf Proof.} Since $V$ is a simple domain in $\C-\Sigma$, there exist 
$|\Sigma|$ nonintersecting 
rays initiating from points of $\Sigma$ and going to infinity,
$\Gamma_1,..,\Gamma_{|\Sigma|}$, such that $\C-\cup_i \Gamma_i$
is a simply connected domain (not containing points from $\Sigma$) and
$V\subset \C-\cup_i \Gamma_i$. Remove now from $\C-\cup_i \Gamma_i$
squares of diameter $2\rho$ centered at points of $\Sigma$, obtaining the 
region $G$. Let $\widehat{V}$ be the 
rectangular domain such that $\widehat{V}\subset U$,
$dist(\widehat{V},\partial U)=\rho$. It is not difficult to represent 
$G\cap \widehat{V}$ as the union of at most $O(|\Sigma|^2)$ quadrilateral domains.
Now slightly perturb the boundary of these domains, so that the resulting
domains $\widehat{V}_1,..,\widehat{V}_M$ are polygonal, contain the unperturbed
domains, and their boundaries constitute a set of $N$ segments, 
$N\leq O(|\Sigma|^2)$, which we call $\g_1,..,\g_N$, such that
$dist(\g_i,Z)\geq \rho/C|Z|$, and $dist(\g_i,\partial U)\geq \rho/2$, 
$i=1,..,N$, where $C$ is some universal constant (it is again not difficult,
 though a bit tedious, to construct such a perturbation).\n
Since the number of zeros of any given branch $\widetilde{f}$ of $f$ on $V$ is
not greater than the sum of the numbers of zeros of corresponding branches 
$\widetilde{f}_1,..,\widetilde{f}_M$ on  $\widehat{V}_1,..,\widehat{V}_M$, 
respectively,
it is also not greater than the sum of the absolute values of variations
of arguments of the branches on $\g_1,..,\g_N$, divided by $2\pi$, by
the argument principle. \ \ \ \ $\Box$\n

{\bf Proof of Theorem 1.9.} Without limiting generality, suppose that 
$K\subset U_1\subset\overline{U_1}\subset U$, where $U_1$ is a 
bounded rectangular domain for which $dist(K,\partial U_1)\geq \rho$
(this makes the proof shorter, as it allows us to
use Lemma 1.11 right away). 

Fix $\l\in F-S$, and let $V\subset K$ be a simple domain in $\C-\Sigma_\l$.
Denote by $\Sigma'_\l,Z'_\l$ the sets $\Sigma_\l\cap U_1, Z_\l\cap U_1$,
respectively.
Note that the cardinality of $Z'_\l$ is uniformly bounded over $\l\in F-S$,
since $F$ and $U_1$ are compact and $b(t,\l)\in \O(U\times W)$.
Denote this bound by $N_1$.

By Lemma 1.11, there exist linear segments $\g_1,..,\g_N$, $N\leq O(N_1^2)$,
all contained in $U_1$,
$dist(\g_i,Z'_\l)\geq \rho/CN_1$, $dist(\g_i,\partial U_1)\geq \rho/2$ 
($C$ being a universal constant which we assume to be larger than $2$),
such that the following holds: the number of zeros of a solution of 1.6)
in $V$ is not greater than the sum of the absolute values of variations
of arguments of some solutions of 1.6) on $\g_1,..,\g_N$.

From Theorem 1.10, this sum of absolute values of variations 
of argument is bounded by
$$
1.7) \ \ \ \ \ 
N\pi(n+1)\left(1+{diam(U_1)\over log(3/2)} \ max_{t\in\cup_i\gamma_i} \ max_j \ 
{|a_j(t,\l)|\over |b(t,\l)|}\right).
$$

Since $dist(\g_i,\partial U_1)\geq \rho/2$, we conclude that not only
for $Z'_\l$ $dist(\g_i,Z'_\l)\geq \rho/CN_1$, but also 
for $Z_\l$ $dist(\g_i,Z_\l)\geq \rho/CN_1$ (recall $C\geq 2$).
Since $N\leq O(N_1^2)$ and $\overline{U_1}\subset U$ is compact, 
from Proposition 1.7 we conclude that 1.7) is uniformly bounded over $\l\in 
F-S$, thereby proving the theorem. \ \ \ \ $\Box$\w

{\bf 2. An asymptotic bound on the number of zeros of solutions of 
linear differential equations which depend regularly on parameters.}\n

The (parameter dependent) Picard-Fuchs equations obtained in section 3 
below for the Abelian integrals are of the form
$$
2.1) \ \ \ \ y^{(n)}+{p_{n-1}(t,\l)\over q(t,\l) }y^{(n-1)}+...+
{p_0(t,\l)\over q(t,\l)}y=0,
$$ 
where $p_i(t,\l)$, $i=0,..,n-1$ and $q(t,\l)$ are polynomials in $t,\l$ with
{\em integral} coefficients. Moreover their degree and height (by which we
mean the maximum of the moduli of their coefficients), as well as 
the order $n$ of 2.1), will admit
asymptotic bounds depending only on the degree of the Hamiltonian $H$
(see section 0). In this section we show how one may obtain a quantitative
version of Theorem 1.9 in this case.\n

So suppose that $p_i(t,\l)$, $i=0,..,n-1$ and $q(t,\l)$ in the equation
2.1) are polynomials with integral coefficients of degree $d$ and 
height $M$.
To avoid problems with notation, it is always assumed below that $d\geq 2$,
$M\geq 2$.

We keep notations of section 1:
$S\subset\C^p$, $p=dim(\l)$, denotes the subset of the parameter
space where $q(\cdot,\l)$ becomes identically zero. For $\l\not\in S$, $Z_\l$
denotes the (now finite) 
set of zeros of $q(\cdot,\l)$, and $\Sigma_\l\subset Z_\l$ 
denotes the set of true singular points of the linear differential 
equation 2.1) for that value of parameter. 

Below, $||\l||$ always denotes the $l_\infty$ norm of $\l$ considered as
a vector in $\R^{2p}\cong\C^p$, and $\rho$ denotes some 
positive number between $0$ and $1$.\n

{\bf Theorem 2.1.} {\it Suppose that the linear differential equation
2.1) depends regularly on $\l$ (in $\C^p$). Then for any $\l\not\in S$,
$||\l||\leq 1$, 
a solution of 2.1) cannot have more than 
$$
n\left({M\over \rho}\right)^{d^{O(p^3)}}
$$
zeros in any simple domain in $\C-\Sigma_\l$, whose distance from $\Sigma_\l$
is at least $\rho$, and which is contained in the unit disc.}\n

Let us remind that according to our convention, 
$O(p^3)$ stands for $c_0+c_1 p^3$, where $c_0,c_1>0$ 
are some {\em universal} (i.e. not dependent on $\l$ or other characteristics 
of the problem) constants.\n

In the proof of Theorem 1.9 the existence of a first-order formula for the set 
$\Gamma$, in the language of ordered rings, 
is rather clear, so we do not write it explicitly.
On the contrary, in the proof of Theorem 2.1, we write an explicit
first-order formula for the maximum of the moduli of the coefficients
of 2.1) on certain subset of $\C\times\C^p$. We then obtain the
value of this maximum from a theorem due to Renegar, regarding
the complexity of quantifier elimination in the first order theory of
the reals. A crucial point here is the integrality of the polynomials
$p_i(t,\l)$, $i=0,..,n-1$ and $q(t,\l)$.\n

As we will be writing first-order formulas, the reader, strictly speaking, 
should recall the
necessary definitions and background of elementary mathematical
logic. We give, though, an intuitive
explanation of what a logical formula defining a semialgebraic set is.

A {\em quantifier free formula} in the language of ordered
rings (this {\em language} is just the following set of symbols: 
$\{0,1,+,\cdot,>\}$), is an expression such as
$$
2.2) \ \ \ \ x+z>0 \ \vee \ xyz^5-5y^2+xz-7>0 \ \wedge \ \neg(xyz=0)   
$$
(strictly speaking, we should have written $1+1+1+1+1$ instead of $5$,
$y\cdot y$ instead of $y^2$, etc.).
The formulas $x+z>0$, $xyz^5-5y^2+xz-7>0$ and $xyz=0$ are sometimes called
the {\em atomic predicates} of the quantifier free formula 2.2), and
below we refer to them as such.
To decipher 2.2), recall that $\wedge$ is the logical 'and', $\vee$ is
the logical 'or', and $\neg$ is the logical 'not'. So the formula 2.2)
in fact defines a subset of $\R^3$
, coordinatized by $(x,y,z)$\footnote{or a subset of any other 
Euclidean space three of whose coordinates are labeled by $x,y,z$.}, 
such that
$$
x+z>0 \ \ or \ \ \left( xyz^5-5y^2+xz-7>0 \ \ and \ \ (not \ \ xyz=0) \ \right)
$$
(recall rules of precedence). In other words, the subset defined
by the formula contains points $(x,y,z)$ for which
either $x+z>0$, {\em or} the following holds: $xyz^5-5y^2+xz-7>0$,
{\em and} $xyz=0$ is {\em not} true.

A (quantifier free) 
formula defining a subset of points $(x,y)\in\R^2$ for which
the truthfulness of $x-y=5$ implies (in the sense of mathematical logic) 
the truthfulness of $y>7$
(i.e. the set $\{(x,y): x-y\neq 5\}\cup \{(x,y): y>7\}$) will be written as
$$
x-y=5 \ \rw \ y>7.
$$
Here $e\rw f$ is a shorthand for $(\neg e)\vee f$.

Let us now write a formula defining the set of $(x,y,z,w)\in\R^4$ such that
there exists $u>0$ for which $xyu+zw+5>0$ and such that for all $v\neq 0$,
$xv+yzw\neq 5$. A possible variant is
$$
(\exists u \ (u>0 \wedge xyu+zw+5>0))\wedge(\forall v \ (v\neq 0 \ \rw \ 
xv+yzw\neq 5)).
$$
Here of course $\exists$ stands for the quantifier 'there exists', 
and $\forall$ stands for the quantifier 
'for all' ($\forall$ may be regarded as a 
shorthand for $\neg\exists\neg$). Thus this formula is not quantifier free.
Note that the quantifiers range over variables which take values in $\R$.
This is why the formula above is called a {\em first order} formula 
(had we allowed quantification over subsets of the reals 
we would get a second order
formula, etc.). 

By the Tarski-Seidenberg principle, for every first order formula 
in the language of ordered rings, there exists a quantifier
free formula defining the same set in a corresponding Euclidean space.
One concludes, therefore, that the set defined by the above 
first-order formula is in fact semialgebraic,
which is not clear {\it a priori}. This is precisely the argument we used
in the proof of Proposition 1.7 to show that the set $\Gamma$ is semialgebraic.
\n
 
Theorem 1.1 from [Re] establishes a quantitative
version of the Tarski-Seidenberg principle, as well as establishing
complexity bounds for the quantifier elimination algorithm which the author 
constructs. 
Cited below is an immediate corollary of that theorem. Let us first set
the framework.\n

Consider the formula
$$
2.3) \ \ \ \ (Q_1 x_1)...(Q_t x_t) \ P(y,x_1,..,x_t),
$$
where $y=(y_1,..,y_l)\in\R^l$, $x_m=(x_{m1},..,x_{mn_i})\in\R^{n_i}$, 
$m=1,..,t$, each $Q_i$ stands either for $\exists$ or $\forall$, and where 
$P(y,x_1,..,x_m)$ denotes the a quantifier free formula with atomic predicates
$$
g_k(y,x_1,..,x_t) \ s_{k} \ 0
$$
$k=1,..,K$, 
where each $g_k(y,x_1,..,x_t)\in\Z[y,x_1,..,x_t]$, $s_{k}$ stand for
either of $=,>,<,\geq,\leq,\neq$ (though strictly speaking, only the symbols
$=,>$ are allowed in the language of ordered rings). 
Suppose now that the sum of the degrees
of the different polynomials among $g_k$, $k=1,..,K$, 
is not greater than $D\geq 2$ (which we call
the {\em total degree} of the quantifier free formula $P(y,x_1,..,x_t)$),
and suppose that the height of the polynomials $g_k$ does not exceed 
$A\geq 2$ for all $k=1,..,K$ (in which case we say that the {\em height}
of $P(y,x_1,..,x_t)$ is not greater than $M$).
\n

{\bf Theorem 2.2 [Re].} {\it The subset of $\R^l$ defined by the formula 2.3),
can be also defined by the quantifier free formula   
$$
2.4) \ \ \ \ \bigvee_{i=1}^I\bigwedge_{j=1}^J \ \ h_{ij}(y) \ s_{ij} \ 0,
$$
such that $h_{ij}(y)\in\Z[y]$ for all $1\leq i\leq I, 1\leq j\leq J$,
the total degree of 2.4) is not greater than 
$$
D^{2^{O(t)}l\Pi_k n_k},
$$
and its height is not greater than}
$$
A^{O(l)+D^{2^{O(t)}\Pi_k n_k}}.
$$

The proof of Theorem 2.1 consists of repeating the proof of Theorem 1.9
together with using Theorem 2.2 to estimate the value of 1.7).\n

{\bf Proof of Theorem 2.1.} We follow here the proof of Theorem 1.9, taking
$U=\C$, taking $K$ to be the origin centered square of diameter $2\sqrt{2}$,
and taking $U_1$ to be the origin centered square of diameter $4\sqrt{2}$.
We prove the theorem for simple domains in $\C-\Sigma_\l$, 
whose distance to $\Sigma_\l$ is at least $\rho$, and which are contained
in the square $K$. In the proof of Theorem 1.9, the estimate for the 
number of zeros in any such domain is given by 1.7)
$$
N\pi(n+1)\left(1+{diam(U_1)\over log(3/2)} \ max_{t\in\cup_i\gamma_i} \ max_j \ 
{|p_j(t,\l)|\over |q(t,\l)|}\right),
$$
where $N\leq O(|Z'_\l|^2)$, $Z'_\l=Z_\l\cap U_1$, 
and $\gamma_i$, $i=1,..,N$ are certain segments 
lying at a distance of at least $\rho/O(|Z'_\l|)$ from $Z_\l$.

Since $|Z'_\l|\leq |Z_\l|\leq d$, to obtain the sought bound for the 
number of zeros from 1.7) it remains to estimate
$$
2.5) \ \ \ \ max_{t\in\cup_i\gamma_i} \ max_j \ 
{|p_j(t,\l)|\over |q(t,\l)|}.
$$
Pick an integer (universal constant) $C>0$, and let $R$ be the nearest integer 
to $1/\rho$, $R\geq 1/\rho$, so that the distance of the segments
$\gamma_1,..,\gamma_N$ to $Z_\l$, which is at least $\rho/O(|Z'_\l|)$, 
is bounded by $1/(CRd)$ from below.
Clearly 2.5) is bounded by the maximum over $j=0,..,n-1$ of
$$
sup_{\l\in\{\l:||\l||\leq 1\}-S} \ sup_{t\in U_1(1/(CRd),\l)} \ \ 
{|p_j(t,\l)|\over |q(t,\l)|},
$$
which by Proposition 1.7 is finite (recall from section 1 that 
$U_1(1/(CRd),\l)$ denotes a subset of points of $U_1$ whose distance
to $Z_\l$ is at least $1/(CRd)$).

We now construct a first order formula for the one point subset of $\R$
whose only point is the value of that maximum. To estimate this value
using Theorem 2.2,
all polynomials entering the formula must have integral coefficients;
this is why $\rho$ was replaced by $1/R$.\n

For a fixed $\l\not\in S$, the set $U_1(1/(CRd),\l)\subset\C\cong\R^2$
is defined by the following formula, which we denote by
$F_1(t_{re},t_{im},\l_{re},\l_{im})$:
$$
\forall(z_{re},z_{im})
\ \ \ ( q_{re}(z_{re},z_{im},\l_{re},\l_{im})=0 \wedge 
\ q_{im}(z_{re},z_{im},\l_{re},\l_{im})=0 
$$
$$
\rw \ CRd((t_{re}-z_{re})^2+
(t_{im}-z_{im})^2)\geq 1 ) \ \ \ \wedge \ \ t_{re}^2\leq 4 \ \wedge \ t_{im}^2\leq 4.
$$ 
The polynomials $q_{re}$ and $q_{im}$ which appear in the formula $F_1$,
are the real and imaginary parts of the polynomial $q$. Observe that
though their degrees are not greater than the degree of $q$, 
their height will be in general
larger than that of $q$ (for example, the height of $t^3$ is $1$, but
each of the real and imaginary parts of 
$$
(t_{re}+\sqrt{-1}t_{im})^3=(t_{re}^3-3t_{re}t_{im}^2)+\sqrt{-1}
(-t_{im}^3+3t_{re}^2t_{im})
$$
has height $3$). To save on space, below we do not write the formulas 
explicitly in terms of the real and the imaginary parts of the variables,
but use, when possible, a shorthand form, 
from which the reconstruction of the formula
is more or less immediate. For example, such a shorthand form of the formula
$F_1$ is (recall that by our conventions 
$||\cdot||$ denotes the $l_\infty$ norm)
$$
\forall z
\  ( q(z,\l)=0  
\rw \ CRd |t-z|^2 \geq 1 ) \ \ \wedge \ \ ||t||\leq 2.
$$
The height of each of $q_{re},q_{im}$ is bounded by 
$2^{2d}(d+1)^{p+1}M$. By Theorem 2.2, there is a quantifier free formula 
$F_2(t,\l)$ equivalent to the formula 
$F_1(t,\l)$ above, its total degree being
not greater than 
$$
d^{2^{O(1)}\cdot (2p+2)\cdot 2}=d^{O(p)},
$$
and its height (recall that 
$d\geq 2$, $M\geq 2$ and $C$ is an integral {\em universal}
constant) is not greater than
$$
max(CRd,2^{2d}(d+1)^{p+1}M)^{O(2p+2)+d^{2^{O(1)}\cdot 2\cdot 2}}
\leq (MR)^{d^{O(1)}(p+1)^2}.
$$
Writing $q(t,\l)=q_0(\l)+q_1(\l)t+..+q_{d}(\l)t^d$, the (quantifier free) 
formula $F_3(\l)$ for the set $\{\l:||\l||\leq 1\}-S$ is given by
$$
||\l||\leq 1 \ \ \wedge \ \ \neg\left(\bigwedge_{i=1}^d q_i(\l)=0\right),
$$
its total degree being not greater than $d$, and its height being not greater
than $2^{2d}(d+1)^{p}M$.

Fix $j\in\{0,..,n-1\}$, and let $v_j$ denote the value of  
$$
2.6) \ \ \ \ sup_{\l\in\{\l:||\l||\leq 1\}-S} \ sup_{t\in U_1(1/(CRd),\l)} \ \ 
{|p_j(t,\l)|\over |q(t,\l)|}.
$$
Since by our assumption, 2.1) depends regularly on parameters,  
such $v_j\in\R$ exists.
We write now a first order formula for the subset $\R-[-v_j,v_j]$ of $\R$:
$$
\forall(t,\l) \ \left(F_2(t,\l)\wedge F_3(t,\l) \ \rw \ 
s^2|q(t,\l)|^2 > |p_j(t,\l)|^2\right).
$$
The total degree of the quantifier free part of this first order formula
is not greater than $d^{O(p)}$, and its height is not greater than 
$(MR)^{d^{O(1)}(p+1)^2}$. Therefore, using 
theorem 2.2 once more, there exists an 
equivalent quantifier
free formula, whose total degree is not
greater than
$$
(d^{O(p)})^{2^{O(1)}\cdot 1\cdot (2p+2)}=d^{O(p^2)},
$$
and its height is not greater than
$$
((MR)^{d^{O(1)}(p+1)^2})^{O(1)+(d^{O(p)})^{2^{O(1)}\cdot (2p+2)}}\leq
((MR)^{d^{O(1)}(p+1)^2})^{d^{O(p^2)}}\leq (MR)^{d^{O(p^3)}}.
$$
We may write now a first order 
formula for the value $v_j$ itself (i.e. for the one point
set $\{v_j\}\subset\R$). Omitting the simple details, we obtain by
Theorem 2.2, a quantifier free formula for the value $v_j$ 
$$
2.4) \ \ \ \ \bigvee_{i=1}^I\bigwedge_{l=1}^L \ \ h_{il}(s) \ s_{il} \ 0,
$$
of total degree and height being bounded again by $d^{O(p^2)}$, 
$(MR)^{d^{O(p^3)}}$.
Since $v_j$ is the only point in the set defined by 2.4), $v_j$ 
must be in fact a root of one of the polynomials $h_{il}$,
which are polynomials (in one variable) with
{\em integral} coefficients. 
Since the degree of this polynomial is bounded by $d^{O(p^2)}$, and its
height is bounded by $(MR)^{d^{O(p^3)}}$, we conclude that $v_j$, for all
$j=0,..,n-1$, is bounded by 
$$
d^{O(p^2)}\cdot (MR)^{d^{O(p^3)}} \ \leq \ (MR)^{d^{O(p^3)}} \ \leq \ 
 \left( M\over \rho\right)^{d^{O(p^3)}}
$$ 
(though $R\geq 1/\rho$, this is not a mistake, 
because of the $O(\cdot)$ notation).
We thus conclude that a bound on the number of zeros may be given by
$$
O(d^2)\cdot (n+1) \cdot O\left( \left( M\over \rho\right)^{d^{O(p^3)}} \right)\ 
= \ n \left( M\over \rho\right)^{d^{O(p^3)}}. \ \ \ \ \Box
$$
\bigskip

{\bf 3. Construction of Picard-Fuchs equations for Abelian integrals.}\n

In this section, we construct a linear 
differential equation, satisfied by the
Abelian integral 
$$
3.1) \ \ \ \ y(t)=\int_{\g(t)} P(x,y)dx+Q(x,y)dy,
$$ 
where $P,Q$ are polynomials,
and $\gamma(t)$ is a continuously varying cycle (see below) 
in the locally trivial 
bundle determined by the mapping $H(x,y):\C^2\rw \C$, where $H(x,y)$ is
also a (generic) 
polynomial. Such linear differential equations are known to exist and
are usually called Picard-Fuchs equations. More precisely, suppose
that the degree of the polynomials $H,P,Q$ is not greater than $d$, and
denote the coefficients of $H(x,y)$ by the tuple $\l$, $dim(\l)=
(d+1)(d+2)/2\leq O(d^2)$.
Then we construct a parameter-dependent 
linear differential equation, with coefficients
being quotients of integral polynomials in $t,\l,\mu$, such that for  
$\l$ not lying in an exceptional codim 1 constructible subset of 
$\C^{dim(\l)}$, 3.1) is a solution of that equation for some value of 
$\mu$.
\footnote{To be precise, this will be true for a constructible subset
of codimension zero in the space of the coefficients of the polynomials $P,Q$.}
Here, $dim(\mu)=(d-1)^2\leq O(d^2)$ as well. 
In section 4 we will show that the linear differential equation we construct
(for any fixed $d$) depends regularly on the parameters $\l,\mu$, 
i.e. is not singularly
perturbed along any arc in the parameter space on which it is defined.
We then show that these equations moreover remain regularly dependent
on parameters also after any rather general algebraic parameter-dependent
change of variable. This property is a direct consequence of the 
algebro-geometric origin of these equations, and is (of course) 
false for general 
regularly dependent on parameter linear differential equations.\n

As in Introduction, 
${\cal H}^d$ denotes the space of polynomials in two variables of
degree not greater than $d$, coordinatized by tuples $\l$ of
polynomials coefficients, so that ${\cal H}^d\cong \C^{dim(\l)}$, $dim(\l)
=(d+1)(d+2)/2$. $H(\l)$ will denote, for $\l\in{\cal H}^d$, the polynomial
whose tuple of coefficients is $\l$. We write either 
$\l\in{\cal H}^d$ or $H\in{\cal H}^d$, depending on the context.\n

We first give the necessary background, explaining the precise meaning 
of the integral 3.1).

It is known, and is not difficult to prove, that for each 
$H\in{\cal H}^d$, there exists a finite set $\Sigma_H\subset\C$,
whose points are called the atypical values of $H$,
such that $H:\C^2-H^{-1}(\Sigma_H)\rw \C-\Sigma_H$ is a 
(smooth) locally trivial bundle. When the projective curve defined by
$H(x,y)=0$, $H$ being a polynomial of degree $d$,
intersects the (complex) line at infinity at precisely $d$
points, we say that the polynomial $H$ is regular at infinity. 
Being $H$ regular at infinity depends only on its highest homogeneous part.
The set of atypical points of a polynomial regular at infinity 
coincides with the set of its critical values. 

Choose $t_0\in\C$ which is not an 
atypical value of $H$, and let $\gamma(t_0)$ 
be a homology class in the first homology group 
of the level curve $\{(x,y):H(x,y)=t_0\}$. It is a basic fact, implied
solely by being $H:\C^2-H^{-1}(\Sigma_H)\rw \C-\Sigma_H$ a smooth
locally trivial bundle, that this homology class admits
a natural continuation along any path in $\C-\Sigma_H$ initiating at $t_0$,
depending only on the homotopy class in $\C-\Sigma_H$ 
of the path joining $t_0$ and the given 
end point. This continuation is given, roughly speaking, by trivializing
the locally trivial bundle  $H:\C^2-H^{-1}(\Sigma_H)\rw \C-\Sigma_H$
at points of the path forming an increasing sequence (with respect to the path
parameter), and transporting a representative of the class from a fiber
to a neighbouring fiber by means of these trivializations. 

Take now a polynomial (or holomorphic) form $P(x,y)dx+Q(x,y)dy$ on $\C^2$, and
let $\gamma(t_0)$ denote, as above, an element in the first homology
group of a point $t_0\in\C-\Sigma_H$. Denote by $\gamma(t)$,
$t\in\C-\Sigma_H$, the continuation of $\g(t_0)$ to $t$ along some path
lying in $\C-\Sigma_H$
($\gamma(t)$, as explained above, is multivalued, depending only on the 
homotopy class of the path connecting $t_0$ to $t$ in $\C-\Sigma_H$).  
The key fact is, that the integral 3.1), called sometimes a complete
Abelian integral, defines in fact an {\em analytic} multivalued function on
$C-\Sigma_H$ ([AGV], Chapter 10). 

In fact, the same considerations apply not only when $t$ varies, but also
when both $t$ and $\l$ vary. Consider the map 
$$
L:\C^2\times{\cal H}^d\rw\C\times{\cal H}^d,
$$
given by $(x,y,\l)\mapsto (H(\l)(x,y),\l)$. It defines a locally trivial
bundle over the open subset $T\subset\C\times{\cal H}^d$, consisting
of pairs $(t,\l)$, for which $H(\l)$ is  
regular at infinity and has degree $d$, and
$t$ is not a critical value of $H(\l)$. Fix $(t_0,\l_0)\in T$
and let $\gamma(t_0,\l_0)$ be some element from the first homology group
of $L^{-1}(t_0,\l_0)\cong\{(x,y):H(\l_0)(x,y)=t_0\}$.
As before, for any holomorphic
form  $P(x,y)dx+Q(x,y)dy$ (or even $P(x,y,\l)dx+Q(x,y,\l)dy$,
where $P,Q$ depend (holomorphically) on $\l$), the integral
$$
3.2) \ \ \ \ \int_{\g(t,\l)} P(x,y)dx+Q(x,y)dy,
$$   
where by $\g(t,\l)$ we denote the multivalued continuation of 
$\gamma(t_0,\l_0)$, will be an analytic multivalued function on $T$.\n

Now fix $H \in{\cal H}^d$. Polynomial forms on $\C^2$ are the same
as regular forms on the quasiprojective variety $\C^2$. The regular 1-forms
on $\C^2$, $\Omega^1(\C^2)$, carry a natural ($H$ dependent) 
structure of $\C[t]$-module, where $\C[t]$ being the ring of 
polynomials in $t$ or, equivalently, regular functions on $\C$ (this structure
is given simply by $t\cdot\o=H(x,y)\cdot\o$; it is useful because it 
commutes with integration: 
$\int_{\g(t)}t\cdot\o=t\cdot \int_{\g(t)}$, where $\cdot$ on the left is
the product operation of the module, and $\cdot$ on the right is the usual
product in $\C$). The submodule {\em generated} by $d\a$, $\a\in\Omega^0(\C^2)$
(i.e. $\a$ is a polynomial in two variables), will be denoted by 
$d\Omega^0(\C^2)=d\Omega^0$. 
The following theorem is essentially contained in [I1].\n

{\bf Theorem 3.1.} {\it Let $H(x,y)$ be a polynomial  
regular at infinity. 
Let $\o\in\Omega^1(\C^2)$. Then the following
are equivalent:\n
i) for all 
continuously varying cycles $\gamma(t)$ in the locally trivial bundle 
determined by $H$, $\int_{\g(t)}\o\equiv 0$,\n 
ii) $\o\in d\Omega^0(\C^2)$.}\n

For our purposes we will need also the following quantitative assertion,
from whose proof one may also extract the proof of Theorem 3.1.\n

{\bf Proposition 3.2.} {\it Let $H(x,y)$ be a polynomial
of degree $d\geq 2$.\n 
i) \ The 1-form $\a(x,y)dx+\b(x,y)dy\in \Omega^1(\C^2)$ 
belongs to
the submodule $d\Omega^0(\C^2)$ 
if and only if there exist $A(x,y),B(x,y)\in\C[x,y]$  
such that 
$$
\a(x,y)dx+\b(x,y)dy=dA(x,y)+B(x,y)dH(x,y).
$$ 
ii) \ Suppose that $H$ is regular at infinity, and suppose that $\a dx+\b dy\in
d\Omega^0(\C^2)$, $deg(\a), \ deg(\b)\leq D$, $D\geq 1$.
Then there exist polynomials $A,B\in\C[x,y]$ of degree not greater than 
$D d^{O(1)}$, such that
$\a dx+\b dy=dA +B dH$.}

{\bf Proof.} Let $\o\in d\Omega^0$. Then, by definition, $\o=\sum_i p_i(t)
\cdot dq_i$,
$p_i(t)\in\C[t]$ and $q_i\in\C[x,y]$ for every i. It suffices therefore to
show that $t^k\cdot dq$ is of the form $dA+BdH$ for all $k$ 
and polynomials $q(x,y)$. But
$$
t^k\cdot dq(x,y)= H(x,y)^k dq(x,y)=d(H(x,y)^k q(x,y))-kH(x,y)^{k-1}q(x,y)dH(x,y).
$$
As for the bound on the degrees of $A,B$, we 
repeat the argument from the proof in [I1] of
Theorem 3.1, but with quantitative estimates.

So suppose that $H\in{\cal H}^d$ is regular at infinity and has degree $d$, 
and suppose $\a(x,y)dx+\b(x,y)dy\in \Omega^1(\C^2)$ 
belongs to
the submodule $d\Omega^0(\C^2)$. 
Without limiting generality, we assume that the leading 
coefficient of $H$ w.r.t. $y$ is $c y^{d}$, $c\neq 0$ (otherwise 
perform a linear change of coordinates to make it such). 
As shown above, $\o=dA+BdH$ for some
polynomials $A,B$. Now $dA+BdH=d(A+BH)-HdB$. For all 
continuously varying cycles $\gamma(t)$ in the locally trivial bundle 
determined by $H$, 
$$
3.3) \ \ \ \ 
\int_{\g(t)}\o \ = \ \int_{\g(t)}(d(A+BH)-HdB) \ = \ \int_{\g(t)}d(A+BH)-t 
\int_{\g(t)}dB \ \equiv \ 0,
$$
(note that this establishes Theorem 3.1 in one direction).\n

Fix now any $x_0\in\C$. The roots
of $H(x_0,y)-t=0$ will be all different for all but finitely many 
(in fact, at most $d-1$) values of $t$. Take $t_0\in\C$ 
for which these roots are all distinct and number the roots as 
$y_1(t_0),..,y_d(t_0)$. Continuing them analytically, we get branches
of analytic multivalued (algebraic) functions $y_1(t),..,y_d(t)$. 

Consider now, for any $t$ which is
not a critical value of $H$,
$$
3.4) \ \ \ \ r(t)=\left( \int_{x_0,y_d(t)}^{x_0,y_1(t)}\o \ \ ... \ \
\int_{x_0,y_d(t)}^{x_0,y_{d-1}(t)}\o \right).
$$  
where the integration is along paths connecting $y_d(t)$ to $y_i(t)$,
$i=1,..,d-1$.
Such paths exist since the affine curve $\{(x,y):H(x,y)-t=0\}$ is 
connected; indeed, this affine curve is smooth since $t$ is not a 
critical value
of $H$, and if it was disconnected, it would imply that the corresponding
projective curve has a singularity on the complex line at infinity -- but
this would contradict the fact that $H$ is regular at infinity. 
Moreover, after making a choice for 
$y_1(t),..,y_d(t)$ (i.e. numbering in some way the roots of $H(x_0,y)-t=0$),
the integrals in 3.4) become well defined, since the integral of $\o$
on the cycle formed by any two paths connecting
$y_d(t)$ to $y_i(t)$, is zero by 3.3).

Consider also
$$
Y(t)=\left( \begin{array}{ccc}
y_1(t)-y_d(t) & \ldots & y_{d-1}(t)-y_d(t) \\
\vdots & \vdots & \vdots \\
y_1^{d-1}(t)-y_d^{d-1}(t) & \ldots & y_{d-1}^{d-1}(t)-y_d^{d-1}(t) 
\end{array} \right),
$$
which is invertible at all points where all roots $y_1(t),..,y_d(t)$ are 
distinct. $r(t)Y(t)^{-1}$ is holomorphic univalued except perhaps for 
the (isolated) points $t$ where
some of the roots $y_1(t),..,y_d(t)$ become equal, or which are critical
values of $H$. Since the growth of $r(t)Y(t)^{-1}$ at these points and
the infinity is polynomial, $r(t)Y(t)^{-1}$ defines in fact a vector
of rational functions in $t$. It is not difficult to see that these rational
functions are in fact polynomials ([Gr], App. C). 
To estimate the degree of these 
polynomials, note that a path from $(x_0,y_i(t))$ to $(x_0,y_d(t))$ 
on $\{H(x,y)-t=0\}$ 
may be constructed as follows.

Consider the projection of $\{H(x,y)-t=0\}\in
\C^2$ to the $x$-plane, and denote the set of critical values 
of this projection
by $Q$. Note that the action 
of the fundamental groop $\pi_1(\C-Q,x_0)$ on the $x_0$-fiber of the projection
is transitive (since the curve is connected and smooth). 
Since the fundamental group is generated by the elementary loops
(i.e. paths which start from $x_0$, go to the vicinity of a point in $Q$,
encircling it and returning to $x_0$ via the same path), and since 
the action is transitive and the fiber consists of $d$ points,
there is a path from $(x_0,y_i(t))$ to $(x_0,y_d(t))$ which projects
to a loop based at $x_0$, composed of traversing at most $d-1$ elementary
loops. We may deform this projection so that it consists of 
at most $6(d-1)$ straight line segments contained in an origin centered disc
which contains all the points of $Q$ 
and lift it to another path, which we denote by $\Gamma$,  
from $(x_0,y_i(t))$ to $(x_0,y_d(t))$ on 
$\{H(x,y)-t=0\}$.
Though the length of $\Gamma$ is not expressible by a semialgebraic
formula, we may estimate the radius of an origin centered ball in 
$\C^2\cong\R^4$, which contains $\Gamma$, and we may also estimate the 
cardinality of the decomposition of $\Gamma$ to connected smooth 
pieces such that the tangent to all the points in a given piece lies in some 
fixed orthant in $\R^4$. Indeed we may write a first-order formula for
the points where $\Gamma$ is either not smooth or the tangent goes from 
one orthant
to another orthant, eliminate the quantifiers by Renegar's theorem, 
and conclude that the 
cardinality of such a decomposition is not greater than $d^{O(1)}$.

Consider now the rate of growth of the maximal norm (in $\C^2\cong\R^4$)   
of the critical points of the projection of $\{(x,y):H(x,y)=t\}$ 
to the $x$-plane.
These correspond to the points of the set $\{(x,y):H(x,y)=t, \ H_y(x,y)=0\}$.
If for $t\rw t_0$, $t_0$ finite, 
one of the points in this set 
tends to $\infty$, this would mean that the projective 
curves defined by $H(x,y)=t$ and $H_y(x,y)=0$ have a common point at infinity,
therefore the projective curves defined by $H_x(x,y)=0$ and $H_y(x,y)=0$ have
a common point at infinity, 
contradicting the assumption that $H$ is regular at infinity. 
So the correspondence $\kappa$ sending $t
\in\C\cong\R^2$ to the maximal modulus
of the critical points of the projection is a semialgebraic 
function which is bounded on compact subsets of its domain. 
Its graph may be given by a quantifier free formula whose 
total degree is independent of $H$ and is bounded by $d^{O(1)}$. 
From Renegar's theorem, for example, one
deduces then that its value at $t\in\C$ is bounded by
$C_H max(2,|t|)^{d^{O(1)}}$, where $C_H>0$ is a constant depending on $H$.
 
By the construction of $\Gamma$, the radius of the origin centered 
ball containing $\Gamma$, is bounded by the maximal distance from
the origin of points of $\{(x,y):H(x,y)=t\}$ projected to the origin
centered disc on the $x$-plane, which has a radius equal to the maximal 
modulus of the critical values of the 
projection plus 1, say. This maximal distance exists for all $t\in\C$, since
the leading coefficient of $H$ w.r.t. $y$ is $cy^{d}$, $c\neq 0$, by
the assumption we made.
Once again, this maximal distance is a semialgebraic function of 
$t\in\C\cong\R^2$ which is bounded on compact sets, and its value
at $t\in\C$ is bounded by $C_H max(2,|t|)^{d^{O(1)}}$,
where $C_H>0$ is a constant depending on $H$.

Therefore the length of $\Gamma$ is not greater than $d^{O(1)}\cdot C_H 
max(2,|t|)^{d^{O(1)}}$, and since the value of the integrals in   
3.4) is bounded by the length of $\Gamma$ times the bound on 
$\sqrt{\|a(x,y)|^2+|\b(x,y)|^2}$ along $\Gamma$, 
where $\a(x,y)dx+\b(x,y)dy=\o$,
we conclude finally that the degree of the entries of 
$r(t)Y(t)^{-1}$ is not greater than $Dd^{O(1)}$.

Writing
$$
r(t)Y(t)^{-1}=(p_1(t),..,p_{d-1}(t))=p(t),
$$
we get
$$
3.5) \ \ \ \ 
\int^{x_0,y_d(t)}_{x_0,y_i(t)}\o-\sum_{l=1}^{d-1}p_l(t)(y_i^l(t)-y_d^l(t)) \
 \equiv \ 0,
$$
$i=1,..,d-1$ (this is possible exactly since $\o\in d\Omega^o(\C^2)$).

Now, 3.5) can be written as 
$$
\int^{x_0,y_d(t)}_{x_0,y_i(t)}\left(\o
-d\left(\sum_{l=1}^{d-1}p_l(H(x,y))y_{i}^l\right)  \right) \
 \equiv \ 0,
$$
which means, putting
$\widehat\o=\o- d\left(\sum_{l=1}^{d-1} 
p_l(H(x,y))y^l\right)$, that the integral
$$
3.6) \ \ \ \ F(u,v)=\int_{(x_0,y_i(H(u,v)))}^{(u,v)}\widehat{\o}
$$
is well defined and holomorphic on $\C^2-H^{-1}(\Sigma_H)$. 
Estimating the growth
of $F(u,v)$ as $u,v$ tend to infinity, we conclude that $F(u,v)$ is in fact
a polynomial of degree $d^{O(1)}$. Indeed, fix $u,v\in\C^2$, and 
consider the path going from $(u,v)$ to 
$(x_0,y_i(H(u,v)))$ on $\{(x,y):H(x,y)-H(u,v)=0\}$, which is 
constructed as the composition of the following two paths. The first path is 
the lift to $\{(x,y):H(x,y)-H(u,v)=0\}$ of the segment lying
 in the $x$-plane which connects $u$ to $x_0$, such that the lifted path
starts at $(u,v)$ and ends at 
$(x_0,y_j(H(u,v))$ for some $j$. Compose now this first path 
with the second path, which connects $(x_0,y_j(H(u,v))$ 
to $(x_0,y_i(H(u,v))$, and which was 
constructed above. We again conclude that the length of the path 
going from $(u,v)$ to 
$(x_0,y_i(H(u,v)))$ on $\{(x,y):H(x,y)-H(u,v)=0\}$ is not greater than 
$d^{O(1)}\cdot max(2,|H(u,v)|)^{d^{O(1)}}$, and it is contained in an
origin centered ball of radius $max(2,|H(u,v)|)^{d^{O(1)}}$. 
When $(u,v)$ tends to a point in $\C^2$, these estimates stay bounded, and 
therefore 3.6) defines a holomorphic function on $\C^2$. When $(u,v)$
tend to infinity, the estimates imply that 3.6) defines in fact a
polynomial of degree $Dd^{O(1)}$ (since the degree of the form $\widehat{\o}$
is not greater than $Dd^{O(1)}$).\n

Writing $F(x,y)$ instead of $F(u,v)$, one has the identity
$$
H_y F_x-H_x F_y=X_H(F)=\widehat{\o}(X_H)=H_y\a-H_x\b,
$$
which implies $H_y(F_x-\a)=H_x(F_y-\b)$. Since $H$ is regular, 
$H_x$ and $H_y$ have no common factor, implying
therefore by unique factorization in $\C[x,y]$, that 
$$
{F_y-\b\over H_y} \ = \ {F_x-\a\over H_x}
$$
is in fact a polynomial, denoted by $B$. It follows that
$\o=dF-BdH$.
Together with the bound $d^{O(1)}$ for the degrees of $F$ and $B$,
this proves the proposition. \ \ \ \ $\Box$\n

Our construction of the Picard-Fuchs equations depends on 
the following result of Gavrilov ([Ga], see 
Proposition 1 and Remark after Lemma 1 in [N]). For us, the degree of 
$P(x,y)dx+Q(x,y)dy\in\Omega^1(\C^2)$ is the maximum of the
degrees of $P(x,y)$ and of $Q(x,y)$.\n

{\bf Theorem 3.3.} {\it Let $H(x,y)$ be a polynomial of degree $d$, 
regular at infinity, 
and let $\widetilde{H}$ denote its highest homogeneous part. 
Then the factor module $\Omega^1(\C^2)/d\Omega^0(\C^2)$
is generated by any set of $(d-1)^2$ forms $\o_1,..,\o_{(d-1)^2}
\in\Omega^1(\C^2)$, for which the polynomials $g_i$, defined by
$$
d\o_i=g_i dx\wedge dy, \ \ \ \ i=1,..,(d-1)^2,
$$
constitute a monomial basis for the complex vector space
$\C[x,y]/<\widetilde{H}_x,\widetilde{H}_y>$. 
Moreover, for any $\o\in\Omega^1(\C^2)$, there exist polynomials 
$a_i(t)\in\C[t]$, $i=1,..,(d-1)^2$, of degree at most $(deg(\o)-deg(\o_i))/d$,
such that $\o=\sum_{i=1}^k a_i(t)\cdot\o_i$.}\n

{\bf Lemma 3.4.} {\it There exists a constructible subset $\G^d\subset
\H^d$ of codimension zero, 
and a set of $(d-1)^2$ monomials $g_1,..,g_{(d-1)^2}
\in\C[x,y]$, their degrees being not greater than $(d-1)^2-1$, such that 
each $H\in\G^d$ is regular at infinity and has degree $d$, 
and for each $H\in\G^d$ $g_1,..,g_{(d-1)^2}$ constitute a basis for
$\C[x,y]/<\widetilde{H}_x,\widetilde{H}_y>$.} 

{\bf Proof.} Observe that for any polynomial $H$ of degree $d$, 
regular at infinity,
the dimension of $\C[x,y]/<\widetilde{H}_x,\widetilde{H}_y>$ is $(d-1)^2$
(cf., for example, [Br] Proposition 2.4).

We assume some familiarity with monomial orderings and leading terms diagrams
(cf. [CLO]). We choose a monomial ordering.
Then, for each $H\in\H^d$, a basis $g_1,..,g_{(d-1)^2}$, of 
$C[x,y]/<\widetilde{H}_x,\widetilde{H}_y>$ may be constructed as follows. 
Let $\Delta_{\widetilde{H}}$ denote
the leading terms diagram of $<\widetilde{H}_x,\widetilde{H}_y>$; 
a basis for $\C[x,y]/<\widetilde{H}_x,\widetilde{H}_y>$ 
may then be taken to consist of  
the $(d-1)^2$ monomials in the complement of $\Delta_{\widetilde{H}}$
([CLO], Pr. 1, pg. 228, Pr. 8, pg. 232). 

Now, there is only a finite number of diagrams with the complement
consisting of $(d-1)^2$ monomials; denote them by $\Delta_1,..,\Delta_l$.
The set ${\cal N}_j\in\H^d$ of polynomials $H$ which are regular at infinity
and have degree $d$, and for which the leading terms diagram of the 
ideal $<\widetilde{H}_x,\widetilde{H}_y>$ is $\Delta_j$, is constructible. 
The sets 
${\cal N}_1,..,{\cal N}_l$ constitute a partition of the set of polynomials
of degree $d$ regular at infinity, which is an open subset of $\H^d$.
Therefore for some $k$, $1\leq k\leq l$, 
${\cal N}_k$ must be of codimension zero. 
We take $\G^d={\cal N}_k$. 

For every $H\in\G^d$, the set of monomials in the complement of
$\Delta_k$, denoted $g_1,..,g_{(d-1)^2}$, will constitute 
a basis of $\C[x,y]/<\widetilde{H}_x,\widetilde{H}_y>$. 
Note the degrees of these monomials cannot be larger than
$(d-1)^2-1$ (because of the structure of leading terms diagrams). 
\ \ \ \ $\Box$\n

In fact, the degrees of the monomials constituting a basis can be taken to be 
not larger than $2(d-1)$ ([N], section 4). 

By Theorem 3.3 and Lemma 3.4, there exist monomial forms $\o_l$, 
$l=1,..,(d-1)^2$, of degree not greater than $(d-1)^2$, 
which generate the $\C[t]$-module $\Omega^1(\C^2)/d\Omega^0(\C^2)$ for any
$H\in\G^d$.\n

We need the following standard fact (10.2.4, pg.284 [AGV]), 
which enables us to write the derivatives
of complete Abelian integrals as complete Abelian integrals again.\n

{\bf Proposition 3.5 (Gelfand-Leray derivative).} {\it Let $H\in\C[x,y]$
and let $\gamma(t)$ denote a continuously varying cycle in the locally
trivial bundle determined by the polynomial $H$. Let $\o$ be a holomorphic
1-form on $\C^2$, and suppose there exists a 1-form $\a$, holomorphic
on $\C^2$, such that $dH\wedge \a=d\o$. Then }
$$
{d\over dt}\int_{\gamma(t)}\o \ = \ \int_{\gamma(t)}\a.
$$
\n

To construct Picard-Fuchs (linear differential) equations, we first 
construct a linear differential system for $(\int_{\g(t)}\o_1,..,
\int_{\g(t)}\o_{(d-1)^2})^T$.\n

Fix $H\in\H^d$, $1\leq l\leq (d-1)^2$. By Proposition 3.5, one can find
polynomials $\a_l,\b_l\in\C[x,y]$ of degree not greater than $d(d-1)$,
so that
$$
3.7) \ \ \ \ {d\over dt}\int_{\gamma(t)} \ (xH_x+yH_y)^2\o_l  \ \ = \ \ 
 \int_{\gamma(t)}\a_l dx+\b_l dy,
$$
holds for any continuously varying cycle $\gamma(t)$ in the locally
trivial bundle determined by $H$.
In fact, we could have taken any polynomial in the ideal $<H_x,H_y>$
instead of $xH_x(x,y)+yH_y(x,y)$. What matters for us, 
however, is that for {\em homogeneous} 
$H$ of degree $d$, $H(x,y)/d=xH_x(x,y)+yH_y(x,y)$.\n

Let now $H\in\G^d$. Note that $H\in\G^d$ has degree $d$ and 
is regular at infinity. By Theorem 3.3 and i) of 
Proposition 3.2 , 
there exist polynomials
$A_l,B_l\in\C[x,y]$, and polynomials $c_{lm}(t)\in\C[t]$, $m=1,..,(d-1)^2$,
so that the integrand on the left hand side of 3.7) can be written as
$$
3.8) \ \ \ \ \left(xH_x(x,y)+yH_y(x,y)\right)^2\o_l 
=  \sum_{m=1}^{(d-1)^2} c_{lm}(H(x,y))\o_m+dA_l(x,y)+B_l(x,y)dH(x,y).
$$
By Theorem 3.3, the degree of $c_{lm}(t)\in\C[t]$ is not greater than 
$(2d+(d-1)^2)/d$, i.e. is not greater than $d$ for $d\geq 2$. 
By ii) of Proposition 3.2, $A_l,B_l$ can be taken to be of degree
not greater than $d^{O(1)}$.
 
Observe that 3.8) defines in fact a linear system over $\C(\l)$, $\l$
being the coefficients of $H\in\H^d$
$$
3.9) \ \ \ \ M(\l)u \ = \ w(\l), 
$$
the unknowns  $u$ being the coefficients of $c_{lm}(t)$ and of $A_l,B_l$,
$1\leq l,m\leq(d-1)^2$. Thus $dim(u)\leq d^{O(1)}$.
The entries of the matrix $M(\l)$ and the vector $w(\l)$
are integral polynomials in $\l$ of degree $d^{O(1)}$ and height not 
greater than $2^{d^{O(1)}}$.
\n

Since 3.9) is solvable for all $\l\in\G^d$, and since $\G^d$ open, 
there exists a constructible set $\G'^d\subset\G^d$ of
codimension zero, and polynomials $p_{lmr}(\l)\in\Z[\l]$,
$l,m=1,..,(d-1)^2$, $r=0,..,d$, and $q(\l)\in\Z[\l]$, such that $\forall \l
\in\G'^d$, $q(\l)\neq 0$, and
$$
3.10) \ \ \ \ c_{lm}(t)=\sum_{r=0}^d {p_{lmr}(\l)\over q(\l)}t^r.
$$
By using Cramer's rule for a suitable subsystem of 3.9), we conclude that
$deg(q(\l))\leq d^{O(1)}$ and its height is not greater than $2^{d^{O(1)}}$,
and the same bounds hold for $p_{lmr}(\l)$. Here we use that 
$dim(\l)\leq O(d^2)$, and the fact that for any two integral polynomials 
$a,b\in\Z[\l]$
$$ 
3.11) \ \ \ \ height(a\cdot b) \ \leq \ (1+min(deg(a),deg(b)))^{dim(\l)}\cdot height(a)\cdot height(b).
$$
Integrating now both sides of 3.8), we get by Theorem 3.1 that
$$
3.12) \ \ \ \  \int_{\g_\l(t)}\left(xH_x+yH_y\right)^2\o_l  
= \int_{\g_\l(t)}\sum_{m=1}^{(d-1)^2} c_{lm}(H)\o_m = 
\sum_{m=1}^{(d-1)^2} c_{lm}(t) \int_{\g_\l(t)}\o_m,
$$
for any continuously varying cycle $\g_\l(t)$ in the locally trivial bundle
determined by $H(\l)$. Using 3.10),
we write 3.12) in the matrix form
$$
3.13) \ \ \ \ \left(\int_{\g_\l(t)} \left(xH_x+yH_y\right)^2\o_m  \right)_m 
\ = \ {K(t,\l)\over q(\l)}\left(\int_{\g_\l(t)}\o_m\right)_m,
$$
where, to save space, $\left(\int_{\g_\l(t)}\o_m\right)_m$ 
denotes the column vector
$\left(\int_{\g_\l(t)}\o_1,..,\int_{\g_\l(t)}\o_{(d-1)^2}\right)^T$, 
and where $q(\l)$
and the entries of $K(t,\l)$ are integral polynomials in $t,\l$ 
of degree and height
bounded by $d^{O(1)}, 2^{d^{O(1)}}$, respectively.\n

Observe that likewise, 
the right hand side of 3.7), may be written using Theorem 3.3,
 Theorem 3.1, and Proposition 3.2, as 
$$
\left(\int_{\g_\l(t)}\a_m dx+\b_m dy\right)_m \ = \ {L(t,\l)\over w(\l)}\left(\int_{\g_\l(t)}\o_m\right)_m,
$$
where, as in 3.13), $w(\l)$ and the entries of $L(t,\l)$ are again 
integral polynomials in $t,\l$ of degree and height
bounded by $d^{O(1)}, 2^{d^{O(1)}}$.\n

The equality 3.7) then implies
$$ 
{d\over dt} \ \ {K(t,\l)\over q(\l)} 
\left(\int_{\g_\l(t)}\o_m\right)_m \ =
 \ {L(t,\l)\over w(\l)} \left(\int_{\g_\l(t)}\o_m\right)_m. 
$$
If the matrix $K(t,\l)$ was invertible over $\C(t,\l)$, we could have
written, for $\l\in\G''^d\subset\G^d$, $\G''^d$ being again a constructible
subset of $\H^d$ of codimension zero, the equality
$$
3.14) \ \ \ \  
\left(\int_{\g_\l(t)}\o_m\right)'_m \ =
 \ K(t,\l)^{-1}\left({L(t,\l)q(\l)\over w(\l)} -{K'(t,\l)} \right)\left(\int_{\g_\l(t)}\o_m\right)_m. 
$$
valid for all continuously varying cycles $\g_\l(t)$ in the locally 
trivial bundle determined by $H(\l)$, $\l\in\G''^d$.

The equality 3.14) would then 
define a linear differential system of dimension
$(d-1)^2$, which would be satisfied,
for any fixed $\l\in\G''^d$, by $\left(\int_{\g_\l(t)}\o_m\right)_m$. 
We write this system as
$$
3.15) \ \ \ \ I' \ = \ {A(t,\l)\over a(t,\l)} I,
$$
where $a(t,\l)$ and the entries of the $(d-1)^2\times(d-1)^2$ matrix
$A(t,\l)$ 
are polynomials in $\Z[t,\l]$. Their degree and height can be again shown,
using 3.11), to be bounded by $d^{O(1)}, 2^{d^{O(1)}}$, respectively.
We note, that since the Abelian integrals may be analytically continued
to all points of the open set $T\subset\C\times\C^{dim(\l)}$, 
3.15) would be in fact
satisfied by $\left(\int_{\g_\l(t)}\o_m\right)_m$ 
for all $\l\in\H^d$ for which $H(\l)$ is 
regular at infinity and has degree $d$, and
for which the matrix $A(t,\l)/a(t,\l)$ is defined 
(as a rational function in $t$).
Thus, the restriction $\l\in\G''^d$ is unnecessary.\n

We now show that the assumption we made is indeed true:\w
{\bf Proposition 3.6.} {\it The matrix $K(t,\l)$ (3.13)) is invertible over
the field $\C(t,\l)$.}

{\bf Proof.} It is for proving this proposition that the integrand
$(xH_x+yH_y)^2\o_m$ in 3.7) was chosen. 
Take any polynomial $H\in\G^d$. By Lemma 3.4 its highest homogeneous part 
$\widetilde{H}$ also belongs to $\G^d$.

Denote the coefficient tuple of $\widetilde{H}$ by $\l_0$.
Since $H\in\G^d$, its highest homogeneous part 
$\widetilde{H}$ is regular at infinity and has degree $d$. Hence, 
any cycle $\g_(t_0,\l_0)$, $t_0\neq 0$, can be naturally continued
to a neighbourhood of $(t_0,\l_0)$.
From Theorem 12.1, pg 317 [AGV], we conclude
that there exists an origin centered $V$ ball in $\C^2$, holomorphic 1-forms 
$\a_1,..,\a_{(d-1)^2}$ on $V$, and cycles $\g_1(t_0,\l_0),..,\g_{(d-1)^2}
(t_0,\l_0)\in H_1(\{H(\l_0)=t_0\}\cap V,\Z)$ for some $t_0\neq 0$ close enough
to $0\in\C$ (i.e. close enough to the critical point of $H(\l_0)$),
for which 
$$
det \ \left(\int_{\g_j(t_0,\l_0)}\a_i\right)_{ij} \ \ \neq \ \ 0.
$$
Instead of the holomorphic forms $\a_1,..,\a_{(d-1)^2}$ 
one may take polynomial forms $\b_1,..,\b_{(d-1)^2}$, 
since on the union of (images of representatives of) 
$\g_1(t_0,\l_0),..,\g_{(d-1)^2}(t_0,\l_0)$, which is a compact subset of
$V$, one may 
approximate $\a_1,..,\a_m$ by polynomial forms to any given accuracy
(true because $V$ is a ball). 

Since $\l_0\in\G^d$, by Theorem 3.3 and Theorem 3.1, 
integral of any polynomial form is a linear combination over $\C[t]$ of
integrals of the polynomial forms $\o_1,..,\o_{(d-1)^2}$ we conclude that
$$
det \ \left(\int_{\g_j(t_0,\l_0)}\o_i\right)_{ij} \ \ \neq \ \ 0.
$$

Now, in a neighbourhood of $(t_0,\l_0)$ 
the matrix $K(t,\l)/q(\l)$ will be equal, as a matrix with {\em meromorphic}
entries, to the holomorphic matrix
$$
\left(\int_{\g_j(t,\l)}\left(xH(\l)_x+yH(\l)_y\right)^2\o_i \right)_{ij}
\cdot \left(\int_{\g_j(t,\l)}\o_i\right)_{ij}^
{-1}.
$$  
But this holomorphic matrix is equal, for $(t_0,\l_0)$ 
to 
$$
\left(\int_{\g_j(t_0,\l_0)}\left(xH(\l_0)_x+yH(\l_0)_y\right)^2\o_i 
\right)_{ij}
\cdot \left(\int_{\g_j(t_0,\l_0)}\o_i\right)_{ij}^
{-1} \ = \ {t_0^2\over d} \cdot Id.
$$ 
precisely since $xH(\l_0)_x+yH(\l_0)_y=H(\l_0)/d$; note that the
right hand side is an invertible matrix (recall $t_0\neq 0$).
But this means that in a neighbourhood of $(t_0,\l_0)$, $K(t,\l)/q(\l)$
is equal (as a meromorphic matrix) to a holomorphic invertible
matrix, implying that $K(t,\l)$ is indeed invertible over $\C(t,\l)$.
\ \ \ \ $\Box$\n

We now show how a linear differential equation for $\int_{\g_\l(t)}\o_m$
may be derived from the system 3.15), $m=1,..,(d-1)^2$.
Taking the derivative of both sides of 3.15), then multiplying by $a(t,\l)$,
and taking into account that $a(t,\l)I'=A(t,\l)I$, we get
$$
a^2(t,\l)I'' \ = \ \left(A^2(t,\l)+a(t,\l)A'(t,\l)-a'(t,\l)A(t,\l)\right)I.
$$
Taking further derivatives and proceeding in the same manner, we get
$$
3.16) \ \ \ \ a^j(t,\l)I^{(j)} \ = \ A_j(t,\l)I,
$$  
where $A_0(t,\l)=Id$, and 
$$
 A_{j+1}(t,\l) \ = \ 
a(t,\l)A'_j(t,\l)+A_j(t,\l)\left(A(t,\l)-ja'(t,\l)Id\right),
$$
$j=0,1,2,.. \ $. Using 3.11), it may be checked that 
$$
3.17) \ \ \ \ deg(A_j)\leq jd^{O(1)}, \ \ \ \ height(A_j)\leq  2^{jd^{O(1)}}
$$
(where $deg(A_j), height(A_j)$ mean the maximal degree and height 
of the entries of $A_j$). Bounds of the same form hold
 for the degree and height of $a(t,\l)^j$.\n

Take now the $m$-th component of $I$, $1\leq m\leq (d-1)^2$. Then for any 
$j=0,1,2,..$, denoting by $\a_{jm}$ the $m$-th row of the matrix $A_j$
in 3.16), we get
$$
3.18) \ \ \ \ a^j(t,\l)I_m^{(j)}=\a_{jm}(t,\l)I.
$$
Clearly there exists $k_m$, $1\leq k_m\leq (d-1)^2$, such that the vectors 
$\a_{0m}(t,\l),..,\a_{k_m m}(t,\l)$ are linearly dependent over $\C(t,\l)$,
but $\a_{0m}(t,\l),..,\a_{(k_m-1)m}(t,\l)$ are not.
Then one may solve uniquely
$$
3.19) \ \ \ \  \a_{k_m m}(t,\l) \ = \ \sum_{l=0}^{k_m-1} 
w_{lm}(t,\l)\a_{lm}(t,\l)
$$
for the rational functions $w_{lm}(t,\l)$, $l=0,..,k-1$.
Multiplying 3.19) by the (column) vector $I$ and using 3.18), 
we get the parameter dependent linear differential equation  
$$
3.20) \ \ \ \ I^{(k_m)}_m-\sum_{l=0}^{k_m-1} 
w_{lm}(t,\l)I^{(l)}_m=0.
$$

{\bf Proposition 3.7.} {\it For all $m=1,..,(d-1)^2$, $k_m=(d-1)^2$.
There exists a constructible set $V_m\subset\H^d$ of codimension zero, 
such that for any $\l\in V_m$, the coefficients of 3.20) are defined, 
and the space of solutions of 3.20) is 
spanned by $\int_{\g_\l(t)}\o_m$.}

{\bf Proof.} Recall that in the proof of Proposition 3.6, we have shown 
that for a homogeneous polynomial $H(\l_0)\in\G^d$, there are 
cycles $\g_1(t_0,\l_0),..,\g_{(d-1)^2}
(t_0,\l_0)\in H_1(\{H(\l_0)=t_0\}\cap V,\Z)$ for some $t_0\neq 0$ close enough
to $0\in\C$ such that 
$$
det \ \left(\int_{\g_j(t_0,\l_0)}\o_i\right)_{ij} \ \ \neq \ \ 0,
$$
and therefore (since $(t_0,\l_0)\in T$), 
$$
3.21) \ \ \ \ det \ \left(\int_{\g_j(t,\l)}\o_i\right)_{ij} \ \ \neq \ \ 0,
$$
in some neighbourhood $U$ of $(t_0,\l_0)$.
By Theorem 3.4 in [AGV], there is $\l_1\in\G^d$ close enough to $\l_0$, 
so that $(t_0,\l_1)\in U$, and 
such that the monodromy 
representation of the fundamental group $\pi_1(\C-\Sigma_{\l_1},t_0)$ on 
the homology group of the fiber $\{(x,y):H(\l_1)(x,y)=t_0\}$ is irreducible.

It is not difficult to show that the homology group 
(with coefficients in $\Z$) 
of any fiber over a point in $T$ is isomorphic to $\Z^{(d-1)^2}$.
Without limiting generality, we may assume that $\g_1(t_0,\l_1),..,
\g_{(d-1)^2}(t_0,\l_1)$ is in fact a basis of 
$H_1(\{(x,y):H(\l_1)(x,y)=t_0\},\Z)$. 
Suppose that $\int_{\g_1(t,\l_1)}\o_m,..,
\int_{\g_{(d-1)^2}(t,\l_1)}\o_m$ are linearly dependent over $\C$ 
(as functions of $t$), so that there exist $\a_1,..,\a_{(d-1)^2}$, 
not all zero, for which
$\sum_{l=1}^{(d-1)^2}\a_l\int_{\g_l(t,\l_1)}\o_m=0$ for all $t$. 
This means that there exists a nonzero element 
$$
\delta=\sum_{l=1}^{(d-1)^2}
\a_l\g_l(t_0,\l_1) \in H_1(\{(x,y):H(\l_1)(x,y)=t_0\},\C)\cong\C^{(d-1)^2},
$$ 
and a nonzero linear form $\g\mapsto\int_\g\o_m$ 
on $H_1(\{(x,y):H(\l_1)(x,y)=t_0\},\C)$,
such that $\delta$ stays in the kernel of this form under all monodromy
transformations.
Therefore the monodromy invariant 
subspace of all elements which remain in the kernel of 
$\g\mapsto\int_\g\o_m$ under all monodromy transformations is nonzero.
Since 
at least one of $\int_{\g_1(t_0,\l_1)}\o_m$,..,
$\int_{\g_{(d-1)^2}(t_0,\l_1)}\o_m$ is nonzero (otherwise the determinant
3.21) would be zero at $(t_0,\l_1)$), this invariant
subspace is nontrivial, and therefore
the monodromy representation 
on $H_1(\{(x,y):H(\l_1)(x,y)=t_0\},\C)$ is reducible, contradicting the
fact cited above. We conclude that $\int_{\g_1(t,\l_1)}\o_m$,..,
$\int_{\g_{(d-1)^2}(t,\l_1)}\o_m$ are linearly independent over $\C$.

We consider now the action of the fundamental group $\pi_1(T,(t_0,\l_1))$
on the homology of the fiber over $(t_0,\l_1)$. 
Choosing a basis for this lattice,  
the monodromy representation acts by linear transformations 
which are represented by integral matrices. Since also the inverse 
transformations are so represented, these matrices
must have determinant equal to $1$ or $-1$, in fact $1$ in our case.
Therefore the Wronskian of 
$\int_{\g_1(t,\l)}\o_m,..,\int_{\g_{(d-1)^2}(t,\l)}\o_m$
is a holomorphic function on $T$, not identically zero. 
By considering its rate of growth as $(t,\l)$ tends to infinity and to points
on $\partial T$ along all complex lines parallel to
the axes of $\C\times\C^{dim(\l)}$, one concludes that the restriction of
the Wronskian to any such complex line is a rational function (in one complex
variable). This implies, however, that the Wronskian is in fact a rational
function in $t,\l$. Thus, there exists a constructible set $V'_m\subset\H^d$ of
codimension zero, such that $\forall \l\in V'_m$ the Wronskian is 
defined and is not identically zero as a function of $t$. We
conclude that there exists a constructible
subset $V_m\subset\H^d$ of codimension zero, 
for which $\int_{\g_1(t,\l)}\o_m,..,\int_{\g_{(d-1)^2}(t,\l)}\o_m$
is a set of linearly independent solutions for 3.20).
 
It follows that $k_m\geq 
(d-1)^2$; since $k_m\leq (d-1)^2$ as well, $k_m=(d-1)^2$.  \ \ \ \ $\Box$\n

Again, since $\int_{\g_\l(t)}\o_m$ continue analytically to all
of $T\subset\C\times\C^{dim(\l)}$, 
it in fact solves 3.20) for all $\l\in\H^d\cong\C^{dim(\l)}$ 
for which $H(\l)$ is regular at infinity and has degree $d$, 
and for which the coefficients of 3.20) are defined.\n

Now, although 
3.20) is a parameter-dependent equation for the integral
$\int_{\g_\l(t)}\o_m$, we are interested in fact in
a linear differential equation for the integral
of an arbitrary 1-form of degree $d$. To obtain such equation
we proceed as follows. 
By Theorem 3.3 and Theorem 3.1 
again, we know that for $\l\in\G^d$, 
the integral of any such form can be written  
as the integral of a linear combination over $\C$ of $\o_1,..,\o_{(d-1)^2}$
(where the coefficients depend on the parameter $\l$).
Let $\mu_1,..,\mu_{(d-1)^2}$ be the coefficients of this linear combination,
which we view again as parameters. Consider now the system 3.15)
augmented by the following equation obtained by differentiating 
both sides of $I_0=\mu_1 I_1+..+\mu_{(d-1)^2} 
I_{(d-1)^2}$, multiplying by $a(t,\l)$, and then simplified 
using 3.15):
$$
 I'_0 \ = \ (\mu_1,..,\mu_{(d-1)^2}){A(t,\l)\over a(t,\l)}I.
$$
We call the new system obtained the {\em augmented system} (for degree $d$).
It depends on the parameters $(\l,\mu)$, and is satisfied by 
(the column vector)
$$
\left(\int_{\g_\l(t)}\sum_{l=1}^{(d-1)^2}\mu_l\o_l,\int_{\g_\l(t)}\o_1,..,\int_{\g_\l(t)}\o_{(d-1)^2}\right)^T  
$$
for any $\l\in\H^d$, for which the system is defined and $H(\l)$
is regular at infinity and has degree $d$.\n

We may now construct a linear differential equation for $I_0$ 
from the augmented system in the same way
by which we constructed the linear differential equations satisfied by the
components of the system 3.15). Its order is {\it a priori} less or equal to 
$(d-1)^2+1$. For generic $\l,\mu$ this equation
has a set of $(d-1)^2$ linearly independent solutions.
Clearly it has also the solution given by a nonzero constant. 
For a generic $\l,\mu$, this constant cannot be
a linear combination of the former solutions by an argument
similar to the argument used in the proof of Proposition 3.7.  
Therefore the order of the equation for $I_0$
is $(d-1)^2+1$. For further reference, 
we write it explicitly as 
$$
3.22) \ \ \ \ I_0^{((d-1)^2+1)}+{p_{(d-1)^2}(t,\l,\mu)\over q(t,\l,\mu)}I_0^{((d-1)^2)}+..+
{p_{0}(t,\l,\mu)\over q(t,\l,\mu)}I_0=0.
$$
$q(t,\l),p_n(t,\l)\in\Z[t,\l,\mu]$, $n=0,..,(d-1)^2$. 
Again, it is satisfied by $\int_{\g_\l(t)}\sum_{l=1}^{(d-1)^2}\mu_l\o_l$
whenever 3.22) is defined for the given values of $\l,\mu$, and
$H(\l)$ is regular at infinity and has degree $d$.\n

{\bf Proposition 3.8.} {\it 
The degree and the height of the polynomials 
$q,p_n\in\Z[t,\l,\mu]$ in 3.22), $n=0,..,(d-1)^2$,
are bounded by $d^{O(1)}, 2^{d^{O(1)}}$, respectively.}\n
{\bf Proof.} A computation by Cramer's rule, taking into consideration 
the bounds for the degree and height of 
$a^j,A_j$, $j\leq (d-1)^2+1$, which are, according to 3.17),
$$
((d-1)^2+1)d^{O(1)}=d^{O(1)}, \ \ \ \ 2^{((d-1)^2+1)d^{O(1)}}=2^{d^{O(1)}}
$$
respectively. \ \ $\Box$\w

{\bf 4. Picard-Fuchs equations depend regularly on
parameters.}\n

Our first aim in this section is to show that the parameter dependent
linear differential equation 3.22) depends regularly on the parameters 
$(\l,\mu)$. In fact, singularly perturbed linear differential equations 
admit parameter dependent solutions with special properties, 
as stated in Lemma 4.1 below. This 
allows to exclude the possibility that such
equations appear as restrictions of Picard-Fuchs equations to
holomorphic arcs in general position in the parameter space. 
We then use Theorem 2.1 (with the 
bounds given in Proposition 3.8) to obtain a bound for 
the number of zeros of an Abelian integral.\n

{\bf Lemma 4.1} {\it Let 
$$
4.1) \ \ \ \ \e^s y^{(n)}+a_{n-1}(t,\e) y^{(n-1)}+...+a_0(t,\e)y=0,
$$
be a linear differential equation such that $s\geq 1$, 
$a_i(t,\e)$ are holomorphic in the 
polydisc $U\times W\subset\C\times\C$, $i=1,..,n$, 
 and for some $k$, $1\leq k\leq n-1$,
$a_k(t,0)\not\equiv 0$ (i.e. 4.1) is, in our terminology, singularly
perturbed). Then for any $t_0\in U$ (with a possible exception of 
a discrete subset of $U$), there exists a parameter dependent
solution of 4.1), $y(t,\e)\in\O(U\times(W-\{0\}))$, for which 
$(y(t_0,\e),..,y^{(n-1)}(t_0,\e))\in\O(W-\{0\})^n$ is a vector with 
constant entries, such that for any other $t_1\in U$ outside of an exceptional,
at most countable, subset $E\subset U$,  
$y(t_1,\e)\in\O(W-\{0\})^n$ has an essential 
singularity at $\e=0$.}

{\bf Proof.} Let $z(t,\e)$ be a parameter dependent solution holomorphic
on $U\times W$ (and not just on $U\times(W-\{0\})$). 
Let $t_0$ be such that $a_k(t_0,0)\neq 0$.

We first show that $z(t,\e)$ is uniquely determined by the $n-1$ functions
$$
y^{(i)}(t_0,\e)\in\O(U\times(W-\{0\})),
$$
$i=0,..,k-1,k+1,..,n-1$. Indeed, if not, 
then there exists such a parameter dependent solution, not 
equal identically to zero, with $z^{(i)}(t_0,\e)\equiv 0$ for
$i=0,..,k-1,k+1,..,n-1$. 

Since $z(t,\e)$ is not identically zero, we then must have that 
$z^{(k)}(t_0,\e)\in\O(U\times(W-\{0\}))$ is not identically zero.
Write $z^{(k)}(t_0,\e)=\e^q
m(\e)$, $q\geq 0, \ m(0)\neq 0$, $m(\e)\in\O(U\times(W-\{0\}))$. 
Suppose $q>0$. 
$z(t,0)$ is the solution of 4.1) at
$\e=0$. Since $a_k(t_0,0)\neq 0$, $z(t,0)$
satisfies a linear differential equation of order not smaller than $k$
and not greater than $n-1$, with initial conditions at $t_0$
being $z^{(i)}(t_0,0)=0 \ \forall i\neq k$,
$z^{(k)}(t_0,0)=(0)^q m(0)=0$, ($q>0$). 
Thus $z(\cdot,0)$ equals
identically to zero. This implies that $z(t,\e)=\e z_1 (t,\e)$,
$z_1(t,\e)$ being another parameter dependent solution satisfying
$z_1^{(i)}(t_0,\e)\equiv 0$ for 
$i=0,..,k-1,k+1,..,n-1$, $z_1^{(k)}(t_0,\e)=\e^{q-1} m(\e)$.
Continuing in the same manner, we finally get the parameter dependent
solution $z_q(t,\e)\in\O(U\times W)$, for which $z_q^{(i)}(t_0,\e)\equiv 0$,
$i=0,..,k-1,k+1,..,n-1$, $z_q^{(k)}(t_0,\e)=m(\e)$. 
Substituting
$z_q(t,\e)$ into the original equation, putting $t=t_0$, we get:
$$
e^s z_q^{(n)}(t_0,\e)+a_{k}(t_0,\e)z_q^{(k)}(t_0,\e)=0, 
$$
which cannot be true, since for $\e=0$ it implies that 
$a_k(t_0,0)\cdot z_q^{(k)}(t_0,0)=a_k(t_0,0)\cdot m(0)\neq 0$ is zero.\n

The set of parameter dependent solutions of 4.1), holomorphic in $U\times W$ 
has the obvious structure of $\O(W)$ module. By the above, this module 
can be identified with a submodule of ${\O}(W)^{n-1}$ ( recall that 
$U\times W$ is connected).

Suppose there exist $n$ such solutions, $y_1(t,\e),..,y_n(t,\e)\in
\O(U\times W)$, linearly independent for at least one value of $\e$ in $W$, 
and identify 
them with the elements $\a_1,..,\a_n\in
{\O}(W)^{n-1}$. 
Then there exist $s_1,..,s_n\in {\O}(W)$, not all
zero, such that $\sum s_i\a_i=0$ (there are such $s_i\in {\M}(W)$, since
$\a_1,..,\a_n$, considered as elements of the ${\M}(W)$ vector space 
${\M}(W)^{n-1}$, are of course linearly dependent; now just multiply by a 
common denominator).
Thus for all $\e\in W$, besides maybe a discrete
set of points where all $s_1,..,s_n$ vanish, 
these $n$ parameter dependent solutions are linearly dependent over $\C$.
Therefore, their (parameter dependent) Wronskian, which is a holomorphic 
function on $U\times W$, vanishes on an open subset of
$U\times W$, implying it in fact vanishes everywhere, and 
there can not exist such a set of parameter dependent 
solutions.

Take now a maximal set of parameter dependent solutions of 4.1),
$y_1(t,\e),..,y_j(t,\e)$, $j<n$,
holomorphic on $U\times W$, which are linearly independent for some parameter
value in $W$. Let us denote this parameter value by $\e_0$, and 
choose a tuple of initial conditions which are linearly independent from
$(y_i(t_0,\e_0),..,y^{(n-1)}_i(t,\e))$, $i=1,..,j$, for example
$(\i_0,..,\i_{n-1})\in\C^{(n-1)}$. Since $j<n$, such a tuple exists.
Let $y(t,\e)$ be a parameter dependent
solution, holomorphic this time only on $U\times(W-\{0\})$, given
by setting all $y^{(l)}(t_0,\e)\in\O(U\times (W-\{0\}))$, $l=0,..,n-1$,
equal identically to  $\i_l$, $l=0,..,n-1$, respectively (it is not
difficult to show, by utilizing the holomorphic dependence on parameter
in $W-\{0\}$ of 4.1), that the solution obtained will be holomorphic on 
$U\times (W-\{0\})$ ).
If $y(t,\e)$ was holomorphic on all of $U\times W$,
we would then get a contradiction to maximality of $y_1(t,\e),..,y_j(t,\e)$
with respect to linear independence. If for some $q>0$ 
$\e^q y(t,\e)$ is holomorphic on $U\times W$, 
it would be another parameter-dependent
solution of 4.1) which would again violate 
the maximality of $y_1(t,\e),..,y_j(t,\e)$
with respect to linear independence. It is now easy to see that $y(t,\e)$
indeed satisfies the conditions of the lemma. \ \ \ \ $\Box$\n

Our argument relies on a 
basic fact about the behaviour of Abelian integrals 
in a neighbourhood of a singularity (a much more precise
information regarding this behaviour is available, but it will not be needed
here). A similar, but not identical, proposition can be found in 
[AGV], Chapter 10.
Below, a {\em meromorphic arc} in $\C^n$ means a (germ of ) 
holomorphic mapping 
from a punctured neighbourhood of $0\in\C$
 into $\C^n$ with at most a pole
at the origin. Thus, a holomorphic arc is also a meromorphic arc, 
but the converse is not necessarily true.\n 

{\bf Proposition 4.2.} {\it Let $(t(\e),\l(\e),\mu(\e))$ be a 
meromorphic arc in the
space $\C\times \C^{dim(\l)}\times\C^{dim(\mu)}$, such that
the arc $(t(\e),\l(\e))$ lies in $T\subset\C\times
\C^{dim(\l)}$ for all $\e\neq 0$, and such that
$$
4.2) \ \ \ \ det \ \left(\int_{\g_j(t(\e),\l(\e))}\o_i\right)_{ij}\ \ \not\equiv 
\ 0
$$
as a function of $t,\e$. Then, for any continuously varying cycle
$\g(t(\e),\l(\e))$ in the locally trivial bundle induced on a 
punctured neighbourhood of $\e=0$ by the locally trivial bundle
over $T$, the Abelian integral
$\int_{\g(t(\e),\l(\e))}\sum_{l=1}^{(d-1)^2} \mu_l(\e) \o_l$ 
can be written in a neighbourhood of the origin 
as a sum of finitely many terms 
$$
4.3) \ \ \ \ \sum_{r,s} \e^{\sigma_r} log^s(\e)h_{rs}(\e),
$$ 
where $h_{rs}(\e)$ are 
holomorphic, and $\s_r$ are certain complex numbers.}

{\bf Proof.} Since $(t(\e),\l(\e))\in T$ $\forall\e\neq 0$, 
$$
\left(\int_{\g_j(t(\e),\l(\e))}\o_i\right)^{-1}_{ij}
$$
is in fact a matrix which is 
analytic multivalued in some punctured neighbourhood of $\e=0$.
Since the monodromy of the derivative (w.r.t. $\e$) of the matrix
$$
\left(\int_{\g_j(t(\e),\l(\e))}\o_i\right)_{ij}
$$
is the same as the monodromy of the matrix itself, and since both 
are analytic (multivalued) in a punctured neighbourhood of $\e=0$,  
$K(\e)$ is in fact holomorphic in this punctured neighbourhood.
The rate of growth of the Abelian integrals as $\e\rw 0$ in any given sector, 
is at most polynomial in $1/|\e|$ (one may estimate the growth by estimating 
the length of transported 
cycles and the diameter of an origin centered ball containing
them, as it was done in the proof of Proposition 3.2; the details are
now different, however, and considerably more tedious - see Extended 
remark 4.3 below). 
Hence $K(\e)$ has at most a pole at $\e=0$. The solution space
of the system
$$
4.4) \ \ \ \ I'_\e \ \  = \ \  K(\e) \ I
$$
is spanned by
$$
\left(\int_{\g_j(t_1,\l(\e))}\o_i\right)_{i}.
$$ 
Since the components of these solutions have at most a polynomial 
rate of growth as $\e\rw 0$ in any given sector, $\e=0$ is a regular 
singular point of 4.4). It is known that solutions of linear
differential systems near a regular singular point are of the form 4.3)
(indeed, the components of such systems satisfy linear differential
equations for which the $\e=0$ is a regular singularity; 
according to Theorems 3.1
and 5.2 in Chapter 4 of [CL], any solution of such equation has the form 4.3)
in the neighbourhood of $\e=0$). Therefore also the integral 
$$
\int_{\g(t(\e),\l(\e))}\sum_{l=1}^{(d-1)^2} \mu_l(\e) \o_l \ = \ 
\sum_{l=1}^{(d-1)^2}\mu_l(\e) \int_{\g(t(\e),\l(\e))}\o_l
$$
can be written in the form 4.3) in a neighbourhood of $\e=0$.
 \ \ \ \ $\Box$\n

{\bf Extended remark 4.3.} In the proof we used the fact  
that the growth of $\int_{\g_j(t_(\e),\l(\e))}\o_i$, as $\e\rw 0$, is at most
polynomial in $1/|\e|$.  Though this fact appears to be well-known, 
it is hard to find a reference for the exact statement; 
we therefore sketch a proof.
Take $\tau\in[0,1]\mapsto \e(\tau)$ be a path in the punctured neighbourhood
of $\e=0$.
A possible way to obtain a continuously varying cycle in the locally trivial 
bundle which is induced over the punctured neighbourhood of zero $E$, starting
from a cycle $\g_0$ in the homology of the fiber 
$\{(x,y):H(\l(\e_0))(x,y)=t(\e_0)
\}$, $\e_0\in E$, is as follows. 

Take a path, say piecewise real analytic, from $\e_0\in E$ to $\e_1\in E$, 
$\tau\in[0,1]\mapsto\e(\tau)\in E$, denoted by $\e(\tau)$.
Take the cycle $\g_0$ in the homology of $\{(x,y):H(\l(\e_0))(x,y)=t(\e_0)
\}$, realized as a real analytic mapping of a circle into
$\{(x,y):H(\l(\e_0))(x,y)=t_1
\}$. Project $\g_0$ on the $x$-axis, obtaining $\delta_0$,
and then construct a continuous deformation 
$\delta(\tau)$, $\delta(0)=\delta_0$, 
so that $\forall \t\in[0,1]$, 
$\delta(\tau)$ never intersects the ($\tau$-dependent) critical values
of the projection of $\{(x,y):H(\l(\e(\tau)))(x,y)=t(\e(\t))
\}$. One can then naturally lift $\delta(\t)$ to $\g(\t)$, a continuously
varying cycle along the path $\e(\t)$ 
in the locally trivial bundle over $E$. One shows that a construction exists
(its precise description being the most tedious part of the proof), 
such that the length of $\delta(\t)$, and consequently of $\g(\t)$, 
is controlled, in a certain precise sense, by the following 
quantities. One is the length of the locus
$L$ traversed by the critical values of the projection to 
the $x$-plane of $\{(x,y):H(\l(\e(\tau)))(x,y)=t(\e(\t))
\}$, as $\t$ varies from $0$ to $1$. The second is
the complexity of the map sending $\t\in(0,1)$ to the collection of the 
critical values of the projection. In the nondegenerate case, this 
complexity is simply the number of self intersections of the locus $L$.

Observe now that to construct a representative for the transport of $\g_0$
from the fiber over $\e_0$ to the fiber over $\e_1$, $|\e_1|\leq |\e_0|$, 
along any path which does 
not wind around $\e=0$, it is sufficient to consider a path composed of
at most four line segments, each of which has a distance of at least
$|\e_1|/\sqrt{2}$ from the origin (since each nonwinding path from $\e_0$
to $\e_1$ is homotopic to such a path). Parameterizing the
line segments in an origin centered disc $D\subset E$, 
by pairs of their endpoints in $D\times D$, 
one then shows, using local finiteness properties of subanalytic 
sets and maps (cf. for example [BM]), that the cardinality of the decomposition of $L$ 
into smooth connected pieces, such that the tangent to each piece 
lies in some fixed quadrant
in $\R^2\cong\C$, is uniformly bounded over all linear segments lying 
in $D$ (we include also segments
which pass through $\e=0$, and for which $L$ may 
have infinite length). 
One also shows that the maximal modulus of the 
critical values of the projection for paths $\e(\t)$ which are line 
segments going
from $\e_0$ to $\e_1$, $|\e_1|\leq |\e_0|$, lying at a distance of at least
$|\e_1|/\sqrt{2}$ from the origin,
is bounded by $C/|\e_1|^\a$, 
for some $C>0,\a>0$. The conclusion is that for such paths, the length of $L$ 
is bounded by $C/|\e_1|^\a$ as well ($C>0,\a>0$ being now some other
constants). The complexity of the locus $L$ for each linear segment lying
in $D$ may be also shown to be uniformly bounded over all such segments,
using similar arguments. One concludes that the length of $\g(\t)$, for each
path going from $\e_0$ to $\e_1$ and which does not wind around the origin,
is bounded by $C/|\e_1|^\a$ for some constants $C>0,\a>0$.

One also shows that the cycle $\g(\t)$ is contained in
an origin centered ball in $\C^2$ of radius at most $C/|\e_1|^\a$ for some,
possibly other, constants $C>0,\a>0$.
Consequently the growth of the integral  
$\int_{\g(t(\e),\l(\e))}\o$, $\o\in\Omega^1(\C^2)$, as $\e\rw 0$ 
(in a given sector), is at most polynomial in $1/|\e|$. (End of Extended 
remark 4.3.)\n

{\bf Theorem 4.4.} {\it The Picard-Fuchs equation 3.22) depends regularly
on the parameters $\l,\mu$ in $\C^{dim(\l)}\times\C^{dim(\mu)} \ 
\cong \ \C^{(d+1)(d+2)/2}\times\C^{(d-1)^2}$.}

{\bf Proof.} 
Suppose not. 
Then there exists a holomorphic arc $(\l(\e),\mu(\e))$
 in the parameter space $\C^{dim(\l)}\times\C^{dim(\mu)}$, such that the 
equation 3.22), restricted to $(\l(\e),\mu(\e))$, (is defined and) 
becomes singularly perturbed. 
Let the set of points of $T\subset\C\times\H^d$
where 4.2) vanishes, be denoted again by $Z$. Let the set of points
of $T\times\C^{dim(\mu)}$, where the Wronskian of
$$
4.5) \ \ \ \ 
1, \ \int_{\g_1(t,\l)}\sum_{l=1}^{(d-1)^2}\mu_l\o_l, \ .. \ , 
\int_{\g_{(d-1)^2}(t,\l)}\sum_{l=1}^{(d-1)^2}\mu_l\o_l
$$
vanishes, be denoted by $Z'$. One concludes from the proofs of Proposition
4.2 and Proposition 3.7 that 
both $Z$ and $Z'$ are intersections with $T$ of proper
algebraic subsets of $\C\times\C^{dim(\l)}$ and $\C\times\C^{dim(\l)}\times
\C^{dim(\mu)}$. 
Therefore, the set $Z''\subset \C\times\C^{dim(\l)}\times
\C^{dim(\mu)}$ defined as the union of the complement of 
$T\times C^{dim(\mu)}$ in $\C\times\C^{dim(\l)}\times
\C^{dim(\mu)}$, with 
$(Z\times\C^{dim(\mu)})\cup Z'$, 
is a closed constructible set of codimension at least 
1 (recall $T$ is constructible and open). 
Let $S\subset\C^{dim(\l)}\times
\C^{dim(\mu)}$ denote the constructible set, such that $(\l,\mu)\in S$
if $(t,\l,\mu)\in Z''$ for all $t\in\C$. 
Since $Z''$ is closed and has codimension at least 1, 
$S$ is closed and has codimension at least 1 as well.

If the arc $(\l(\e),\mu(\e))$ intersects $S$ at infinitely many points $\e_n$,
$\e_n\rw 0$, this means (since $S$ is constructible and closed) that 
the arc lies on $S$ (at least when the parameter 
of the arc is restricted to a sufficiently small disc). Since $S$ has 
codimension at least 1, 
using Lemma 1.6 it is possible to find a nearby holomorphic arc not 
lying on $S$, 
such that 3.22), restricted to this new arc, (is defined and) stays
singularly perturbed. Restricting the arc parameter to lie in a sufficiently
small disc, this arc may intersects $S$ only at $\e=0$.
So we may assume that if the arc 
$(\l(\e),\mu(\e))$ intersects $S$, it intersects it only at $\e=0$.
If the arc did not intersect $S$, 3.22) 
restricted to  $(\l(\e),\mu(\e))$ would have a linearly independent 
set of $(d-1)^2+1$ 
parameter dependent solutions, holomorphic on a polydisc 
(recall that $S$ is closed).
From the proof of Lemma 4.1 it would then follow that 3.22), 
restricted to  $(\l(\e),\mu(\e))$, is not singularly perturbed. 
Therefore $(\l(0),\mu(0))\in S$.

Note that being 3.22) defined on the original arc (which may lie on $S$),
does not imply (at least {\em a priori}) that
the Wronskian of 4.5) is not identically zero when restricted to that arc.
This is the reason that the set $Z'$ appears in the definition of the 
exceptional set $S$. The reason that the set $Z$ appears in this definition,
is to let us use Proposition 4.2 in the argument below. 

Take now a polydisc $U\times W\subset\C\times\C$, such that the
coefficients of 3.22) are holomorphic on $U\times (W-\{0\})$, and
on which the equation 3.22), restricted to $(\l(\e),\mu(\e))$, takes the form
4.1). According to Lemma 4.1, there exists then a parameter dependent 
solution $y(t,\e)\in\O(U\times (W-\{0\}))$ of this equation, 
such that for some $t_0, t_1\in U$, $y^{(j)}(t_0,\e)\in\O(W-\{0\})$, 
$j=0,..,(d-1)^2$, are constants, and $y(t_1,\e)\in\O(W-\{0\})$ 
has an essential singularity at $\e=0$. 
We may assume that $W$, $t_0$ and $t_1$ are such 
that the arcs $(t_0,\l(\e),\mu(\e))$ and $(t_1,\l(\e),\mu(\e))$, where 
$\e\in W$, intersect $Z''$ only at $\e=0$ (indeed, the arc $(\l(\e),\mu(\e))$ 
intersects $S$ only at $\e=0$, and Lemma 4.1 shows that
we may choose $t_0,t_1\in U$ almost arbitrarily). We may then write, for all 
$\e\in W-\{0\}$
(since for $\e\neq 0$ the integrals are defined and the Wronskian is not 
identically zero), 
$$
4.6) \ \ \ \  y(t_0,\e) = \ c_0(\e)\cdot 1 + \sum_{i=1}^{(d-1)^2} 
c_i(\e)\int_{\g_i(t_0,\l(\e))}\sum_{l=1}^{(d-1)^2}\mu_l(\e)\o_l.
$$
Fixing a basis of the homology of the fiber over $(t_0,\l(\e_0))$ for some 
$\e_0\in W$, $\e_0\neq 0$, 
$c_i(\e)$ are then certain branches of analytic multivalued functions
on $W-\{0\}$.
To compute $c_i(\e)$, $i=0,..,(d-1)^2$, fix $t=t_0$. 
$(c_1(\e),..,c_{(d-1)^2}(\e))^T$ is then given by the product of
the inverse of the Wronskian {\em matrix} of 4.5),
and of $(y(t_0,\e),..,y^{((d-1)^2)}(t_0,\e))^T$. 
Note that the latter does not depend on $\e$. 
Since the arcs $(t_0,\l(\e),\mu(\e))$ and $(t_1,\l(\e),\mu(\e))$ 
intersect $Z''$ only at $\e=0$, 
we conclude from Cramer's rule and Proposition 4.2, that $c_i(\e)$ can be
written as ratio of finite sums of the form 
$\sum_{i,j} \ \e^{\sigma_i} log^j(\e)h_{ij}(\e).$
Now, since the arc 
$(\l(\e),\mu(\e))$ intersects $S$ only at $\e=0$, 
the set of points $(t,\e)\in U\times W$, mapped by 
$(t,\l(\e),\mu(\e))$ to points of $Z''$, is of (complex) codimension 1 
(consequently not separating $U\times W$).
By analytic continuation, 4.6) holds therefore also if we replace $t_0$ by  
$t_1$, implying
$$
 \ \ \ \  y(t_1,\e) = \ c_0(\e)\cdot 1 + \sum_{i=1}^{(d-1)^2} 
c_i(\e)\int_{\g_i(t_1,\l(\e))}\sum_{l=1}^{(d-1)^2}\mu_l(\e)\o_l.
$$ 
Using Proposition 4.2 again, we then conclude that $y(t_1,\e)$ can be 
written as a ratio of sums 
$\sum_{i,j} \  \e^{\sigma_i} log^j(\e)h_{ij}(\e)$ as well.
But this is impossible, since $y(t_1,\e)$ has an essential singularity at 
$\e=0$.
Indeed, it is not difficult to show that a ratio of sums 
$\sum_{i,j} \  \e^{\sigma_i} log^j(\e)h_{ij}(\e)$
can never be equal to a holomorphic function in a punctured neighbourhhood
of the origin with an essential singularity at the origin. 

We conclude that 3.22) must depend regularly on the parameters $\l,\mu$ 
in $\C^{dim(\l)}\times\C^{dim(\mu)}$.
\ \ \ \ $\Box$\n

Recall now that for $H\in\H^d$, $\Sigma_H$ denotes the set of its 
atypical points (which are just the critical values for the 
polynomials $H$ we consider below).
\n

{\bf Corollary 4.5.} {\it Let $\l\in\H^d$, $d\geq 2$, 
be such that $H(\l)$ is regular at infinity, and suppose that
$||\l||\leq 1$. Then for any polynomials 
$P,Q\in\H^d$, the integral $\int_{\g(t)}P(x,y)dx+Q(x,y)dy$ can have not
more than
$$
 \left({2\over \rho}\right)^{2^{d^{O(1)}}}
$$
zeros in any simple domain in $\C-\Sigma_H$, whose distance from $\Sigma_H$
is not smaller than $\rho$, $0<\rho<1$, and which is contained in the 
unit disc. Here $\g(t)$ is any continuously varying cycle in the locally
trivial bundle determined by $H(\l)$.} 

{\bf Proof.} 
The coefficients of 3.22) are ratios of integral polynomials 
of degree and height not greater than $d^{O(1)},2^{d^{O(1)}}$, respectively.
The order of 3.22) is $(d-1)^2+1$. For $\l,\mu$ in a certain constructible 
set $V\subset\H^d$ of codimension zero, 3.22) is defined and its solutions 
contain all integrals of the form 
$\int_{\g(t)}\sum_{l=1}^{(d-1)^2}\mu_l\o_l$. Let $(\l,\mu)\in V$ and
suppose $||\l||\leq 1$, $||\mu||\leq 1$.  
The bound for their number
of zeros in any simple domain of $\C-\Sigma_H$, whose distance from 
$\Sigma_H$ is not smaller than $\rho$, and which is contained in the unit disc,
is then given by Theorem 2.1 as (note that the dimension of the parameter
space is $dim(\l)+dim(\mu)\leq (d+1)(d+2)/2+(d-1)^2\leq O(d^2)$)
$$
((d-1)^2+1)\left({2^{d^{O(1)}}\over \rho}\right)^{d^{O(O(d^2)^3)}}
\leq \left({2\over \rho}\right)^{2^{d^{O(1)}}}.
$$
One may omit now the restriction $||\mu||\leq 1$, since the zeros 
of $\int_{\g(t)}\sum_{l=1}^{(d-1)^2}\mu_l\o_l$ are the same as the zeros
of $\int_{\g(t)}\sum_{l=1}^{(d-1)^2}(\mu_l/||\mu||)\o_l$.
Since any integral $\int_{\g(t)}\sum_{l=1}^{(d-1)^2}\mu_l\o_l$ has a natural
analytic continuation to $T\times\C^{dim(\mu)}$, we immediately conclude 
(using Rouche theorem, for example), that for any $\mu$ and $||\l||\leq 1$, 
such that $H(\l)$ is regular at infinity and has degree $d$, the same 
bound for the number of zeros holds. 

Observe that for $\l\in\G^d$ any integral 
$\int_{\g(t)}P(x,y)dx+Q(x,y)dy$ can be written as 
$\int_{\g(t)}\sum_i \mu_i\o_i$ for some $\mu\in\C^{(d-1)^2}$ (Theorem 3.3), 
implying the bound for such integrals as well.
Now, the set $\G^d$ is of codimension zero, and any integral 
$\int_{\g(t)}\sum_{l=1}^{(d-1)^2}Pdx+Qdy$ continues
analytically to $T\times\C^{dim(\mu)}$. This implies that the bound holds 
for the number of zeros of any
integral $\int_{\g(t)}P(x,y)dx+Q(x,y)dy$, where $\g(t)$ is a 
continuously varying cycle in the locally trivial bundle determined by $H(\l)$,
$||\l||\leq 1$, for which $H(\l)$ is regular at infinity and has degree $d$.
\ \ \ \ $\Box$\n

Note that for general 
linear differential equations which depend regularly on parameters, 
making a parameter dependent
algebraic change of variable in general will not produce an equation
which depends regularly on parameters. For example,
$$
{dy\over dt} \ - \ {y \over t(t-\e)} \ = \ 0
$$
becomes, after putting $t=\e\t$:
$$
{dy\over d\t} \ - \ {y\over \e\t(\t-1)} \ = \ 0.
$$
On the contrary, Picard-Fuchs equations stay regularly dependent
on parameters also after a parameter dependent algebraic change of variable.
This is of course the consequence of their algebro-geometric origin, or,
more to the point, the consequence of the fact that the rate of
growth of Abelian integrals along holomorphic curves in the new coordinates,
stays at most polynomial. We now prove a precise claim of this sort.\n

By a rational map from $\C^m$ to $\C^n$ we mean a map whose 
domain is a dense subset of $\C^m$ and whose components are rational
functions.

Let $\b(\k):\C^{dim(\k)}\rw \C^{dim(\l)}\times\C^{dim(\mu)}$ be a dominant
rational map (i.e. a rational map which maps its domain to a dense subset
of the target space). We may write the linear differential equation 3.22)
in the new coordinate and parameters, obtaining (the derivative is with 
respect to $\t$)  
$$
4.7) \ \ \ \ I_0^{((d-1)^2+1)}+{a_{(d-1)^2}(\t,\k)\over b(\t,\k)}I_0^{((d-1)^2)}+..+
{a_{0}(\t,\k)\over b(\t,\k)}I_0=0,
$$
where $b(t,\k),a_i(\t,\k)\in\C[\t,\k]$, $i=0,..,(d-1)^2$.\n

{\bf Proposition 4.6.} {\it The linear differential equation 4.7)
depends regularly on $\k$ in $\C^{dim(\k)}$.}

{\bf Proof.} 
We proceed as in the proof of Theorem 4.4.
Suppose 4.7) does not depend regularly on parameters. 
Then there exists a holomorphic arc $\k(\e)$
 in the parameter space $\C^{dim(\k)}$, such that the 
equation 4.7), restricted to $\k(\e)$, is (defined and) singularly
perturbed. 

Let the set $K\subset \C\times\C^{dim(\k)}$ be the closed constructible subset
obtained as the union of the points where the rational map 
$(\t,\k) \ \mapsto \ (\a(\t,\k), \b(\k))$
is undefined, and of the preimage, under this map, of the set $Z''\subset
\C\times\C^{dim(\l)}\times\C^{dim(\mu)}$ which was constructed in
the proof of Theorem 4.4. Since the derivative of $\a(\t,\k)$ 
is not identically zero and $\b(\k)$ is dominant, $K$ is of codimension
at least 1. Let $S\subset\C^{dim(\k)}$ denote the constructible set, 
such that $\k\in S$ if $(\t,\k)\in K$ for all 
$\t\in\C$. It is closed and of
codimension at least 1 as well.
As in the proof of Theorem 4.4, it may be assumed that 
the arc $\k(\e)$ intersects $S$ only at $\e=0$.
 
Take now a polydisc $U\times W\subset\C\times\C$, such that the
coefficients of 4.7) are holomorphic on $U\times (W-\{0\})$, and
on which the equation 4.7), restricted to $\k(\e)$, takes the form
4.1). Again, 
Lemma 4.1 implies that there exists a parameter dependent 
solution $y(\t,\e)\in\O(U\times (W-\{0\}))$ of this equation, 
such that for some $\t_0, \t_1\in U$, $y^{(j)}(\t_0,\e)\in\O(W-\{0\})$, 
$j=0,..,(d-1)^2$, are constants, and $y(\t_1,\e)\in\O(W-\{0\})$ 
has an essential singularity at $\e=0$. 

We proceed now as in the proof of Theorem 4.4. 
We may again assume that $W$, $\t_0$ and $\t_1$ are such 
that the arcs $(\t_0,\k(\e))$ and $(\t_1,\k(\e))$, where 
$\e\in W$, intersect $K$ only at $\e=0$. We then write, for all 
$\e\in W-\{0\}$
(since for $\e\neq 0$ the integrals are defined and the Wronskian is not 
identically zero), 
$$
4.8) \ \ \ \ y(\t_0,\e) \ = \ \ c_0(\e)\cdot 1 + \sum_{j=1}^{(d-1)^2} 
c_j(\e)\int_{\g_j(\a(\t_0,\k(\e)),
\b_\l\circ\k(\e))
}\sum_{l=1}^{(d-1)^2}(\b_\mu\circ\k(\e))_l\o_l.
$$
The coefficients $c_i(\e)$, $i=0,..,(d-1)^2$, are given, as before, 
by the product of the inverse of the Wronskian matrix of
$$
1, \ \int_{\g_1(\a(\t_0,\k(\e)),
\b_\l\circ\k(\e))}\sum_{l=1}^{(d-1)^2}(\b_\mu\circ\k(\e))_l\o_l, \ .., \ 
\int_{\g_{(d-1)^2}(\a(\t_0,\k(\e)),
\b_\l\circ\k(\e))} \sum_{l=1}^{(d-1)^2}(\b_\mu\circ\k(\e))_l\o_l,
$$
and of $(y(\t_0,\e),..,y^{((d-1)^2)}(\t_0,\e))^T$, the latter being a constant
vector in $\C^{(d-1)^2+1}$. From Proposition 4.2 
we conclude that also in this case $c_i(\e)$ can be
written as ratio of finite sums of the form 
$\sum_{i,j} \ \e^{\sigma_i} log^j(\e)h_{ij}(\e)$.

Since the arcs $(\t_0,\k(\e))$ and $(\t_1,\k(\e))$, where 
$\e\in W$, intersect $K$ only at $\e=0$, 
4.8) holds also if we replace $\t_0$ by  
$\t_1$, implying
$$
\ \ \ \ y(\t_1,\e) \ = \ \ c_0(\e)\cdot 1 + \sum_{j=1}^{(d-1)^2} 
c_j(\e)\int_{\g_j(\a(\t_1,\k(\e)),
\b_\l\circ\k(\e))
}\sum_{l=1}^{(d-1)^2}(\b_\mu\circ\k(\e))_l\o_l.
$$
Again, one concludes that $y(\t_1,\e)\in\O(W-\{0\})$ can be 
written as a ratio of sums 
$\sum_{i,j} \  \e^{\sigma_i} log^j(\e)h_{ij}(\e)$. This is 
impossible, since $y(\t_1,\e)$ has an essential singularity at 
$\e=0$. 

We conclude that 4.7) depends regularly on $\k$ in $\C^{dim(\k)}$.
\ \ \ \ $\Box$\n

We use this result only to derive Corollary 4.7 below, though more
general statements can be made. Together with Corollary 4.5, it
proves Theorem 0.2.\n

{\bf Corollary 4.7.} {\it The condition $||\l||\leq 1$ in Corollary 4.5
may be removed.}

{\bf Proof.} Indeed, let $\a(\t,\k)=\t$ and let $\b(\k_i)=\l_i$ for all $i$,
$2\leq i\leq dim(\l)+dim(\mu)$, $\b(\l_1)=1/\k_1$. Clearly $\b$ is then a
dominating rational map, and $\a(\t,\k)$ has a nonzero derivative w.r.t. $\t$.
By Proposition 4.6, the equation 4.7) depends then regularly on $\k$.
This implies
(since the degree and the height of 4.7) are the same then as of 3.22))
that Corollary 4.5 holds also for all $\l$, for which
$||(1/|\l_1|,\l_2,..,\l_{dim(\l)})||\leq 1$. Considering not just the 
transformation above, but all 
transformations of the form $\l=(\k_1^{s_1},..,\k_{dim(\l)}^{s_{dim(\l)}})$,
for different choices of $s_i\in\{-1,+1\}$, $i=1,..,dim(\l)$, we conclude
that the condition $||\l||\leq 1$ may be indeed removed. \ \ \ \ $\Box$\w

\begin{center}{\bf REFERENCES}\end{center}
{\bf [Aea]}. \ \ \ \  Arnol'd V.I. et al, Some unsolved problems in
the theory of differential equations and mathematical physics. {\it Uspekhi Mat. Nauk} 44
(1989), no. 4(268),191--202; translation in {\it Russian
Math. Surveys} 44 (1989), no. 4, 157--171.\n  
{\bf [AGV]}. \ \ \ \ Arnol'd, V. I.; Gusein-Zade, S. M.;
Varchenko, A. N. {\it Singularities of differentiable
maps. Vol. II. Monodromy and asymptotics of integrals.} 
Monographs in Mathematics,
83. Birkhauser Boston, Inc., Boston, MA, 1988.\n 
{\bf [BM]}. \ \ \ \  Bierstone, E.; Milman, P. \  Semianalytic and 
subanalytic sets. {\it Inst. Hautes \'Etudes Sci. Publ. Math.} 
No. 67 (1988), 5--42.\n 
{\bf [Br]}. \ \ \ \ Broughton, S. A. Milnor numbers and the topology of 
polynomial hypersurfaces. {\it Invent. Math.} 92 (1988), no. 2, 217--241.\n
{\bf [CL]}. \ \ \ \ 
Coddington, E. A.; Levinson, N. {\em Theory of ordinary differential equations.} McGraw-Hill Book Company, Inc., New York-Toronto-London, 1955.\n
{\bf [CLO]}. \ \ \ \ Cox, D.; Little, J.; O'Shea, D. {\it Ideals,
varieties, and algorithms. An introduction to computational algebraic
geometry and commutative algebra.} 
Second edition. Undergraduate Texts in Mathematics. Springer-Verlag, New York, 1997.\n
{\bf [Ga]}. \ \ \ \  Gavrilov, L. Petrov modules and zeros of
Abelian integrals. {\it Bull. Sci. Math.} 122 (1998), no. 8, 571--584.\n 
{\bf [GI]}. \ \ \ \ Glutsuk A.; Il'yashenko Y. 
An estimate of the number of zeros of Abelian integrals
for special Hamiltonians of arbitrary degree. 
{\it Preprint} arXiv:math.DS/0112156v1.\n
{\bf [Gr]}. \ \ \ \ Grigoriev, A. {\it Ph.D. thesis}, the Weizmann Institute
of Science, December 2001.\n 
{\bf [I1]}. \ \ \ \ Il'yashenko, Yu. Appearance of limit cycles in perturbation of
the equation ${dy\over dz}=-{R_z\over R_w}$, where $R(z,w)$ is a polynomial.
{\it USSR Math. Sbornik}, 78 (1969), 360-373.\n
{\bf [I2]}. \ \ \ \ Il'yashenko, Yu. Centennial history of Hilbert's 16th problem. {\it Bull. Amer. Math. Soc. (N.S.)} 39 (2002), no. 3, 301--354.\n 
{\bf [IY]}. \  Il'yashenko, Y.; Yakovenko, S. Double exponential
estimate for the number of zeros of complete Abelian integrals and rational envelopes of linear
ordinary differential equations with an irreducible monodromy group.
{\it Invent. Math.} 121 (1995), no. 3, 613--650.\n
{\bf [Kh]}. \ \ \ \ A. Khovanskii.  Real analytic manifolds with the
property of finiteness, and complex abelian integrals, {\it
Funktsional. Anal. i Prilozhen.} 18 (1984), no. 2, 40-50.\n
{\bf [N]}. \ \ \ \ Novikov, D. Modules of the abelian integrals and the 
Picard-Fuchs systems. {\it Nonlinearity} 15 (2002), no. 5, 1435--1444.\n
{\bf [NY1]}. \ \ \ \ Novikov, D.; Yakovenko, S. Tangential Hilbert
problem for perturbations of hyperelliptic Hamiltonian systems. 
{\it Electron. Res. Announc. Amer. Math. Soc.} 5 (1999), 55--65
(electronic).\n 
{\bf [NY2]}. \ \ \ \ Novikov, D.; Yakovenko, S. 
Trajectories of polynomial vector fields and ascending chains of 
polynomial ideals. {\it Ann. Inst. Fourier (Grenoble)} 
49 (1999), no. 2, 563--609.\n
{\bf [NY3]}.\ \ \ \ 
Novikov, D.; Yakovenko, S. Redundant Picard-Fuchs system for abelian integrals. {\em J. Differential Equations} 177 (2001), no. 2, 267--306.\n
{\bf [R]}. \ \ \ \ Roussarie, R. {\it Bifurcation of planar vector fields and Hilbert's sixteenth problem.} 
Progress in Mathematics, 164. Birkhäuser Verlag, Basel, 1998.\n 
{\bf [Re]}. \ \ \ \ Renegar, J. On the computational complexity and
geometry of the first-order theory of the reals. III. Quantifier
elimination. {\it J. Symbolic
Comput.} 13 (1992), no. 3, 329--352.\n
{\bf [Var]}. \ \ \ \ A.N. Varchenko. Estimation of the number of zeros of an
abelian integral depending on a parameter, and limit cycles, {\it
Funktsional. Anal. i Prilozhen.} 18 (1984), no.2, 14-25.\n 
{\bf [Y]}. \ \ \ \ Yakovenko S. On functions and curves defined by ordinary differential equations
, to appear in {\it Proceedings of the Arnoldfest} (Ed. by Bierstone,
Khesin, Khovanskii, Marsden), Fields Institute Communications, 1998.\n

\end{document}